\documentclass{article}

\usepackage{graphicx}
\usepackage{amsmath}
\usepackage{amssymb}
\usepackage{bm}
\usepackage{float}
\usepackage{multicol}
\usepackage{tabularx}
\usepackage{xcolor}
\usepackage{comment}

\newcommand{\R}{\mathbb{R}}
\newcommand{\mH}{\mathrm{H}}

\title{Least Squares with Equality constraints Extreme Learning Machines for the resolution of PDEs}
\author{Davide Elia De Falco, Enrico Schiassi, Francesco Calabrò}
\date{June 2025}

\begin{document}

\maketitle

\tableofcontents

\section{Abstract}
In this paper, we investigate the use of single hidden-layer neural networks as a family of ansatz functions for the resolution of partial differential equations (PDEs). In particular, we train the network via Extreme Learning Machines (ELMs) on the residual of the equation collocated on -eventually randomly chosen- points. Because the approximation is done directly in the formulation, such a method falls into the framework of Physically Informed Neural Networks (PINNs) and has been named PIELM. Since its first introduction, the method has been refined variously, and one successful variant is the Extreme Theory of Functional Connections (XTFC). However, XTFC strongly takes advantage of the description of the domain as a tensor product. Our aim is to extend XTFC to domains with general shapes. The novelty of the procedure proposed in the present paper is related to the treatment of boundary conditions via constrained imposition, so that our method is named Least Squares with Equality constraints ELM (LSE-ELM). An in-depth analysis and comparison with the cited methods is performed, again with the analysis of the convergence of the method in various scenarios. We show the efficiency of the procedure both in terms of computational cost and in terms of overall accuracy. 

\section{Introduction}
A numerical method for the resolution of Partial Differential Equations (PDEs) typically relies on describing the numerical solution as a linear combination of (piecewise-) polynomial functions. The use of diverse ansatz functions to solve PDEs has a long and varied history, characterized by fluctuating opinions and mixed results. Without aiming to be exhaustive, the selection of different function spaces is usually driven by the desire to avoid grid generation, mitigate the Runge phenomenon, or impose expected or known properties of the problem's solution on the function space.

In this paper, we investigate the use of single hidden layer neural networks as ansatz functions. In particular, we train the network via Extreme Learning Machines (ELMs) on the residual of the equation collocated on (eventually randomly chosen) points. The training of ELMs degrees of freedom is done by the resolution of a linear problem, the neurons being fixed in a random way and not tuned. Randomness in the generation of the neurons gives the resolution benefits of some overparameterization and that the resolution is done via Least Squares (LS). When Neural Networks (NN) are used to approximate the solution of a PDE, and the approximation is made directly as a choice of test functions, the method is called physic-based 
and the structures are called PINNs \cite{cuomo}.  In recent years, many researchers have used ELMs for the resolution of Differential Equations in such PINNs framework. In the next subsection, we present an extensive review of the available research in such an area: the use of ELMs in PINNs. What we can conclude is that despite the increasing number of papers, some of the fundamental understandings are still lacking. Also, most methods are well suited only for some specific problems.  

{The starting points for the study of the present paper are PIELM, as named in \cite{Dwivedi2020}, and their alternative named XTFC in \cite{Schiassi2021}. PIELM, as already pointed out in \cite{Calabro2021}, is a collocation method that is eventually solved in the LS sense. XTFC adds to this method the use of Functional Connections to treat the boundary conditions analytically. Thanks to this, many problems can be efficiently solved, especially in the 1D setting (ordinary differential problems, two points boundary problems, optimal control) or with tensor product domains but equivalent approaches for general domains are still lacking. 
In the present paper, we propose a different approach to overcome this limitation} and give insights concerning the notion of convergence in such a setting. The last is an important theoretical and practical question that was already addressed in \cite{defalco2024} in the case of interpolation, whereas generally the focus on the use of NN for the resolution of PDEs only focuses on accuracy. The question well arises in the construction with ELM because no iteration is needed for the training as the overall training problem is linear. In the present paper, the convergence is discussed when the training set is increased (more collocation points); when the number of degrees of freedom is increased (more neurons in the hidden layer); and finally when both are increased keeping fixed the ratio between the two.

This gives that we have the possibility to obtain accurate solutions that in many PINN strategies are not obtained, due to the fixed number of degrees of freedom and the suboptimality of nonconvex nonlinear optimizations. The structure of ELMs, which in our personal vision bridge spectral methods and "pure" PINNs, allows an easy definition of spectral convergence. In our numerical tests, we can also see nice convergence and high accuracy in challenging problems, both linear and nonlinear, both stationary and time-dependent.  The method is also tested in cases with minor regularity, and in this case, convergence is not guaranteed.

Moreover, it is well known that the treatment of boundary conditions is crucial to obtain a versatile and accurate method for PDEs. Boundary treatment in collocation and in PINNs methods has a vast literature, so that also in this case we present in a subsection the literature review. 

Now, let us introduce some preliminary notation.
    In the rest of the paper we will focus on problems of the type
    \begin{equation}\label{eq:Pb1}
        \begin{cases}
            \mathcal D [u](\bm{x}) = f(\bm{x}) \qquad \bm{x} \in \Omega \\         
            \mathcal B [u](\bm{x}) = g(\bm{x}) \qquad \bm{x} \in \Gamma
        \end{cases}   
    \end{equation}
    where {$\mathcal D$ is a differential operator (linear or nonlinear), $\mathcal B$ is a linear operator,} $\Omega \subset \mathbb R^d$ is an arbitrary domain with $d=2,3$ and $\Gamma\subset \mathbb R^{d-1}$; usually $ \Gamma = \partial \Omega$. {In the rest, we will not distinguish between spatial and temporal variables, since the whole treatment will be carried out in the same way. Namely, we consider space-time discretizations and treat the initial as part of the boundary conditions, so that in these cases $\Gamma$ wil not include points at final time.}

The discrete solution that we look for is denoted by $\Tilde{u}\left(\bm{x}\right)$, and is described as linear combinations of $L$ - in our case fixed- functions $\sigma_j$ combined via the weights $\beta_j$: $\Tilde{u}\left(\bm{x}\right)=  \sum_{j=1}^L \beta_j \sigma_j(\bm{x})$. 
\\
The parameters $\bm{\beta}\in \mathbb R^L$ are the unknowns, determined via some exactness condition imposed on a set of $N$ points, where in general we have $N_I$ of them inside $\Omega$ and $N_B$ on  $\Gamma$.

Both the determination of the functions $\sigma_j$ and the discussion on the final discrete analogous of problem \eqref{eq:Pb1} are discussed in the next section.

The remainder of the paper is organized as follows. In the next two subsections we review the literature focusing on the main two features explored in the present paper: the treatment of boundary conditions and the use of ELMs for the resolution of PDEs. In Section \ref{sect:PRELIM} we introduce and analyze our method. 
In Section \ref{sect:COMPARISON} we use our method, PIELM and XTFC to solve a simple problem where all three can be applied. Then, we make an extensive comparison both on the convergence properties and on the spectral behavior of the discrete counterpart of the PDE. 
In Section \ref{sect:TEST} we test our procedure in general-shaped domains and in the case of challenging problems, selected among the most well-known test examples mainly taken from \cite{CaloAnisotropicDiffusion,collection2Dproblems}.  
Finally, in Section \ref{sect:Conclusions} we draw some conclusions. 

\subsection{Treatment of the boundary conditions in collocation methods: Review of the literature}
Incorporating boundary conditions in collocation methods has a relevant literature due to the interest in such methods in different settings. 
In the Neural Network application to the resolution of PDEs via PINNs, the issue of the incorporation of boundary conditions has been widely considered. Being the training of a network an iterative procedure, many methods focus on the possibility of setting such conditions iteratively. \\
The first paper dealing with NN training for the resolution of differential problems is \cite{lagaris1998}
where the solution is sought as the sum of two contributions, one of the two introduced so as to impose exactly the boundary condition. 
After this, most of the existing PINN approaches enforce the essential (Dirichlet) boundary conditions by means of additional penalization terms that contribute to the loss function, weighted via some scaling parameter. 
In order to tune these scaling parameters during the training, such techniques require involved implementations and, in most cases, lead to intrusive methods since the optimizer has to be modified.
We refer to \cite{berrone2023enforcing} for a review of the imposition of essential boundary conditions in PINNs, and to \cite{deng2024physical} for some general boundary terms.
All these impositions of boundary conditions seem to bring some new error in the overall procedure due to the additional terms included in the optimization. Also, in some cases, the imposition of such conditions needs to be solved iteratively, and this increases the computational cost. 

In the case of ELM training, the construction of the parameters that define the solution is done directly, so that we need to impose also the boundary conditions once, and then the case is much similar to the other more traditional collocation methods such as mesh-free, spectral, or isogeometric. 
The first issue deals with the imposition of essential boundary conditions (Dirichelet), for which we refer to \cite{wagner2000application} for a nice introduction to the analysis in a mesh-free environment. Essential boundary conditions are naturally imposed only considering the behavior on the boundary points, so that in collocation methods that include the evaluation at boundary points seem more reliable to use. Despite this, in \cite{wagner2000application} the authors suggest that the only application of the conditions, named in the paper "traditional collocation", cannot restore optimal properties, while the transformation of the bases made to satisfy the conditions is to be preferred. The second issue is in the case of the imposition of natural (Neumann) boundary conditions in collocation schemes. In this case, many efforts have been made when the structure of the space is well known, such as spectral polynomial methods \cite{Wang2010} or isogeometric collocation methods \cite{de2015isogeometric}. This is because
the approximation of the boundary terms also takes into account the behavior at internal points.
\\
In summary, many alternatives have been introduced, but all can be divided using the classification proposed in \cite{berrone2023enforcing} and we can distinguish soft, hard, and weak imposition of the boundary conditions. \\
\begin{itemize}
    \item Soft imposition is obtained by adding a penalization term in the loss function in such a way that the optimizer will converge to a solution that satisfies the boundary conditions as well, this is the case of \cite{Dwivedi2020}. The main advantage of this approach is that only the loss function needs to be modified, resulting in a very easy and efficient implementation. 
    \item Hard imposition is obtained by modifying the structure of the neural network in order to satisfy the boundary conditions by construction. A typical example is the addition of a non trainable layer at the end of the network as considered in \cite{Schiassi2021}
    \item Weak imposition is obtained by enlarging the space of ansatz functions and imposing a penalization term in the weak formulation of the residual as shown in \cite{berrone2023enforcing}.
\end{itemize}

\subsection{The use of ELMs for the resolution of PDEs: Review of the literature}
The first paper that uses a method that relies on ELM for the resolution of PDEs is \cite{Sun2018}. In this paper the ELM are not considered as the ansatz functions.
The first articles that introduce the use of ELM networks as ansatz functions in the resolution of linear PDEs are \cite{Calabro2021,Dwivedi2020}, dealing with both linear problems and the soft imposition of Dirichelet boundary conditions. In particular, \cite{Calabro2021} focuses on elliptic problems with boundary layers, while \cite{Dwivedi2020} results in 2D general-shape domains and discusses the accuracy and versatility of the procedure. Then, in the same year Schiassi et al. introduce the XTFC in \cite{Schiassi2021} and test accuracy when the boundary condition imposed hardly. XTFC are then used in various PDE applications, see \cite{Schiassi2021, Schiassi2022, Schiassi2022b, Wang2024}. 

Starting from these results, ELMs have been used in different applications. 
Linear elasticity problems are considered in \cite{Yan2022, Zhu2024}. The Stokes equations are solved in \cite{Chen2022}.
Higher order PDEs problems have been considered in \cite{Dwivedi2020b, XiAnPRE, sun2025two}.
Nonlinear PDEs have been considered firstly in \cite{Dong2021, Dwivedi2021, Fabiani2021} and more recently in \cite{Harder2024,Neufeld2024, Thiruthummal2024}. The moving boundaries for the Stefan problem are considered in \cite{ren2025physics}.
Inverse problems are considered in \cite{Dong2023}.

Various variants of ELMs are also tested for the resolution of linear and nonlinear PDEs. As an example, in \cite{Calabro2023} the authors explore the possibility of using ELMs for spatial discretizations and standard multistep methods for time derivatives in parabolic problems. Other attempts include the use of IELM (Incremental ELM) \cite{Yang2020}; HELM (Hierarchical ELM) \cite{DongPRE}; MLELM (Multi-Layer ELM) \cite{Dong2022,Ni2023}; KELM (Kernel ELM) \cite{Li2023}; Piecewise ELM \cite{LiangPRE}; Chebychev ELM \cite{Liu2021}.

In most of the previously cited papers, the focus is on the accuracy of the numerical solution, and only in few cases convergence is considered, see \cite{Calabro2021,Calabro2023,Dong2021,Dong2022,Dong2023,Fabiani2021,Ni2023,Yan2022}.

It is worth to notice that before using ELM in PDEs, ordinary differential problems are tackled in \cite{Yang2018} with soft imposition of the initial condition and used also in \cite{Panghal2020}. First use of hard imposition of initial condition can be found in \cite{Schiassi2021b} following the XTFC setting and has been used in \cite{Daryakenari2023, DeFlorio2022, DeFlorio2024, DeFlorio2024b, Laghi2023, Schiassi2024}. Recently some results also on generalization errors have been derived in the same setting, see \cite{Schiassi_etalPRE}. Such estimates are also used for the resolution of optimal control problems in \cite{Schiassi2024}.

Moreover, also two-point boundary problems have been considered with ELMs \cite{DAmbrosio2021, DAmbrosio2024, DeFlorio2021} and recently also fractional ODEs in \cite{Kumar2024, Kumar2024b, Sivalingam2024}.



\section{Introduction to LSE-ELM and relations with PIELM and XTFC}\label{sect:PRELIM}
In this section we introduce the proposed method. For clarity and to fix the notation, we start with an introduction to ELMs in Section \ref{Sect:ELMs}, and then move to the treatment of the operators in Section \ref{Sect:coll}.  The treatment of the boundary condition is in Section \ref{Sect:boundary}, where we also describe the PIELM and the XTFC to facilitate a rapid comparison, and our proposal is detailed in Section \ref{Sect:Constrained}. 

\subsection{ELM Neural Networks}\label{Sect:ELMs}
In this paper we focus on linear scalar equations in two dimensions, while reporting theoretical results, we emphasize that such results could be reported in a more general framework. A Single Layer Feedforward Neural Network (SLFN) with $L$ hidden neurons is determined by the input weights $\{\bm{w}_j\in\R^d\}_{1\,\leq\, j\,\leq L}$, the hidden neurons' biases $\{b_j\in\R\}_{1\,\leq\, j\,\leq L}$, the output weights $\{\beta_j\in\R\}_{1\,\leq\, j\,\leq L}$, and an activation function $\sigma:\R\to\R$. In what follows, we restrict our attention to the case in which $\sigma$ is the hyperbolic tangent and the $\bm{w}_j$ and $b_j$ are sampled independently from the uniform distribution on $[-M,M]$, where the parameter $M$ controls how steep are the resulting basis functions. \\
In order to fix the notation for the training of the SLFN, let us consider the following approximation problem, see \cite{defalco2024}.
 Given $N$ samples $(\bm{x}_i,y_i) \in \R^d \times \R$, the SLFN aims to approximate the data $y_1,\ldots,y_N$ with weighted sums; in formula:
        \begin{equation}
            \sum_{j=1}^L \beta_j\, \sigma(\bm{w}_j^{\mathrm{T}}\bm{x}_i+b_j) \approx y_i\,,\quad 
i = 1\,\ldots,{N}\,,
        \end{equation}
        where $\bm{\beta}=[\beta_1,\ldots,\beta_L]^T\in\R^L$ is referred to as output weights vector.\\  
ELM is a training algorithm for SLFNs where the input weights $\{\bm{w}_j\}_{{1 \leq j\leq L}}$ and the biases of hidden neurons $\{b_j\}_{{1 \leq j\leq L}}$ are randomly selected from a continuous distribution and remain fixed, the optimization is performed exclusively on the output weights vector. Since the training set is also fixed, the entries of the matrix 
        \begin{equation}\label{eq:Hij}
            \mH_{ij} = \sigma(\bm{w}_j^{\mathrm{T}}\bm{x}_i + b_j)\,.
        \end{equation}
are fixed as well and $\bm{\beta}$ can be efficiently computed as the solution of a linear optimization problem of the form
        \begin{equation}
            \min_{\bm{\beta} \,\in\, \R^L} \| \bm{\mH \beta} - \bm{y} \|\,,
        \end{equation}
        with $\|\cdot\|$ denoting the Euclidean norm, and
        where $\bm{\mH}\in\R^{n\times L}$ is defined as the output matrix of the hidden layer, a typical approach is using the Moore-Penrose pseudoinverse:
        \begin{equation}    \label{MP_solution}
            \bm{\beta} = \bm{\mH}^\dagger \bm{y}\,,
        \end{equation}
     Notice that this is not the sole choice to obtain the solution to the minimization problem. In fact, being the problem of Least-Square type, one can consider other valid resolution strategies such as the ones related to the resolution of the normal equations or orthogonal triangularization, however we found that the Moore-Penrose pseudoinverse (specifically the formulation based on the SVD decomposition of the matrix $\bm{\mH}$) has the best balance between accuracy on training data and generalization error on new data. 

\subsection{Discretization of the operator: collocation}\label{Sect:coll}
    We consider as collocation points $\{\bm{x}_i\}_{{1 \leq i\leq N_I}}$ in $\Omega$ , and get a discrete version of problem \eqref{eq:Pb1}:
    \begin{equation}
        \label{eq:Pb1_collocation}
        \mathcal D [u](\bm{x}_i) = f(\bm{x}_i) \qquad 1 \leq i\leq N_I
    \end{equation}

    As an approximation of the solution $u: \Omega \rightarrow \mathbb R$ we consider an ELM, so that $\sigma_j(\bm{x})=\sigma(\bm{w}_j\bm{x} + b_j)$ and:
    \begin{equation}
        \Tilde{u}\left(\bm{x}\right) = \sum_{j=1}^L \beta_j \sigma_j(\bm{x})
    \end{equation}
    In the case of a linear operator, problem \eqref{eq:Pb1_collocation} can be written as:
    \begin{equation}
        \sum_{j=1}^L \beta_j \mathcal D [\sigma_j](\bm{x}_i) = f(\bm{x}_i) \qquad 1 \leq i\leq N_I
    \end{equation}
    At this point, the operator acts on $L$ fixed functions evaluated on $N_I$ collocation points, and we can introduce $\bm{A}$ and $\bm{f}$ such that:
    \begin{equation}
        A_{ij} = \mathcal D [\sigma_j](\bm{x}_i) \qquad f_i = f(\bm{x}_i)
    \end{equation}
    We can finally write problem \eqref{eq:Pb1_collocation} as:
    \begin{equation}
        \label{eq:Pb1_discrete}
        \bm{A\beta} = \bm f
    \end{equation}
    which can be interpreted as an ELM whose activation functions have been modified by the differential operator.\\
  In the case of a nonlinear operator the training is a nonlinear optimization with the following residual:
    \begin{equation}
        \mathcal R_i(\bm \beta) = \mathcal D[\Tilde{u}](x_i) - f(x_i) \qquad 1\leq i \leq N_I
    \end{equation}

\subsection{Imposition of boundary conditions}\label{Sect:boundary}
In this section, we first review the method for the hard imposition of boundary conditions. In fact, a typical goal of linear models in this context is to separate the interpolation problems in $\Omega$ and in $\Gamma$, that is to say:
    \begin{enumerate}
        \item Find any function $u_B$ that satisfies the boundary conditions
        \item Fix $u_B$ and consider the modified problem
            \begin{equation}
                \begin{cases}
                    \mathcal D [u_I+u_B](\bm{x}) = f(\bm{x}) \qquad &\bm{x} \in \Omega \\         
                    \mathcal B [u_I](\bm{x}) = 0 \qquad &\bm{x} \in \partial \Omega
                \end{cases}
            \end{equation}
        \item Choose $u_I$ as a linear combination of functions that satisfy homogeneous boundary conditions.
        \item Set the final approximation to be $\Tilde{u} = u_I + u_B$
    \end{enumerate}

    When $\Omega = [x_0,x_f] \times [y_0,y_f]$, Extreme Theory of Functional Connections (XTFC) is an effective way to deal with linear boundary conditions. As shown in \cite{Schiassi2021} we can introduce a linear transfinite interpolator $\mathcal M$ (such as Coons' Patches) and define:
    \begin{align}
        u_B &= \mathcal M[g](\bm x) \\
        u_I &= \sum_{j=1}^L \beta_j \Bigl( \sigma_j(\bm x) - \mathcal M[\sigma_j](\bm x)\Bigr) = \sum_{j=1}^L \beta_j\Tilde{\sigma_j}(\bm x)
    \end{align}
    where $u_B$ is just an interpolation of the prescribed data from $\Gamma$ to $\Omega$, and in order to satisfy the homogeneous boundary conditions we can subtract to any $\sigma_j$ the interpolation of its own values on $\Gamma$. \\
    Then $\forall \bm{x} \in \Gamma,\, \forall \bm{\beta} \in \mathbb R^L$:
    \begin{align*}
        \mathcal B[\Tilde{u}](\bm{x}) &=  \mathcal B[u_I](\bm{x}) + \mathcal B[u_B](\bm{x}) \\
        &= 0 + \mathcal B[u_B](\bm{x}) \\
        &= g(\bm{x})
    \end{align*}
   {For linear operators} we can introduce $\bm A \in \R^{N_I \times L}$ and $\bm f \in \R^{N_I}$ such that:
    \begin{equation}
        A_{ij} = \mathcal D[\Tilde{\sigma_j}](\bm x_i) 
        \qquad
        f_i = f(\bm x_i) - \mathcal D[u_B](\bm x_i)
    \end{equation}
    Therefore we have the following collocation problem in the interior of $\Omega$:
    \begin{equation}
        \bm{\beta}^* = \arg\min \frac{1}{2}\|\bm{A\beta} - \bm f\|_2^2
    \end{equation}
 
    As already mentioned, the case of a nonlinear operator leads to a nonlinear optimization problem:
    \begin{align*}
        \mathcal R_i(\bm \beta) &= \mathcal D[u_I+u_B](\bm x_i) - f(\bm x_i) \qquad 1 \leq i \leq N_I \\
        \bm \beta^* &= \arg\min \frac{1}{2}\|\bm {\mathcal R}(\bm \beta)\|_2^2 \ .
    \end{align*} 
    
    It is important to note these two properties of the XTFC:
    \begin{enumerate}
        \item The loss function depends on $\bm{\beta}$ only through $u_I$ and the dependence is still linear, and the model can still be trained in the ELM framework. 
        \item The additional time required for the modifications is negligible, since analytical expressions for Coons' Patches are available on tensor product domains.
    \end{enumerate}
  The previous strategy can be generalized to convex polygons \cite{Salvi2014}, and using different interpolation strategies it is possible to work on arbitrary polygons.
    Strategies that allow to work on arbitrary domains \cite{Dyken2009} exist but unfortunately, analytical expressions for the interpolants are not available and their evaluations lead to quadratures or numerical optimization problems that are harder and require more computations than the original problem. On top of all of these complications, everything scales up with the dimensions of the problem \cite{Floater2008}.

We are now interested in the introduction to soft imposition of boundary conditions, that is probably the most common approach in the PINN literature, as it asks only to perform collocation on the boundary. 
This has also been used for ELMs as Physics-Informed ELM (PIELM) in \cite{Dwivedi2020}:  the collocation points $\{\bm{x}_k\}_{{1 \leq k\leq N_B}}$ on $\Gamma$ are added to the training set and we define $\bm C \in \R^{N_B \times L}$ and $\bm g \in \R^{N_B}$ such that:
    \begin{equation}
        C_{kj} = \mathcal B [\sigma_j](\bm{x}_k) \qquad g_k = g(\bm{x}_k)
        \qquad
        1 \leq k \leq N_B, \, 1 \leq j \leq L
    \end{equation}
The collocation for linear and nonlinear problems is as follows:
    \begin{align}
        \bm{\beta}^* &= \arg\min \Bigl(\frac{1}{2}\|\bm{A\beta} - \bm f\|_2^2 + \frac{\lambda}{2}\|\bm{C\beta} - \bm g\|_2^2 \Bigr) \\
        \bm{\beta}^* &= \arg\min \Bigl(\frac{1}{2}\|\bm{\mathcal R}(\bm \beta)\|_2^2 + \frac{\lambda}{2}\|\bm{C\beta} - \bm g\|_2^2 \Bigr) 
    \end{align}
    In general, $\lambda$ is a penalization term and needs to be tuned, but in most applications the accuracy and generalization capability of ELMs is enough even if $\lambda$ is not optimal and many works PIELM is employed with $\lambda=1$.
    However, if very accurate results are needed, tuning $\lambda$ is important and as it grows the least squares problem gets more unstable, and this could be detrimental for the generalization on unseen data. \\
If the problem of the training is seen as LS with equality contrains, the strategy above is the popular weighting method \cite{Stewart1997}.

\subsection{Constrained imposition of boundary conditions: the proposed method}\label{Sect:Constrained}
    In order to reformulate the problem and impose all kind of boundary conditions in all kind of domains, we notice that the two equations can be collocated on the two set of nodes (as in the case of soft imposition), but one of the two can be taken as a constraint for the other:
    \begin{equation}
        \min_{\bm{C\beta} = \bm g}\|\bm{A\beta}-\bm{f}\|
    \end{equation}
    If the problem of the training is seen as LS with equality contraints, this is equivalent to a weighting method with $\lambda \rightarrow \infty$. Also, this has a strict relation with the hard imposition, being the result of the method a modification of the feature functions, as we will point out.



\begin{figure}[H]
        \label{fig:projected_outputs}
        \centering
        \includegraphics[width=0.45\textwidth]{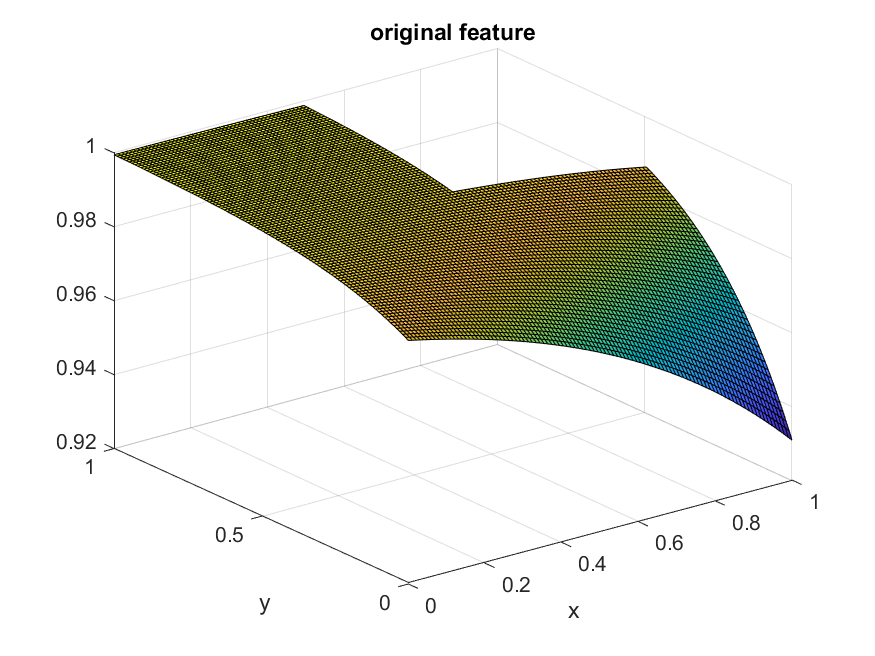}
        \includegraphics[width=0.45\textwidth]{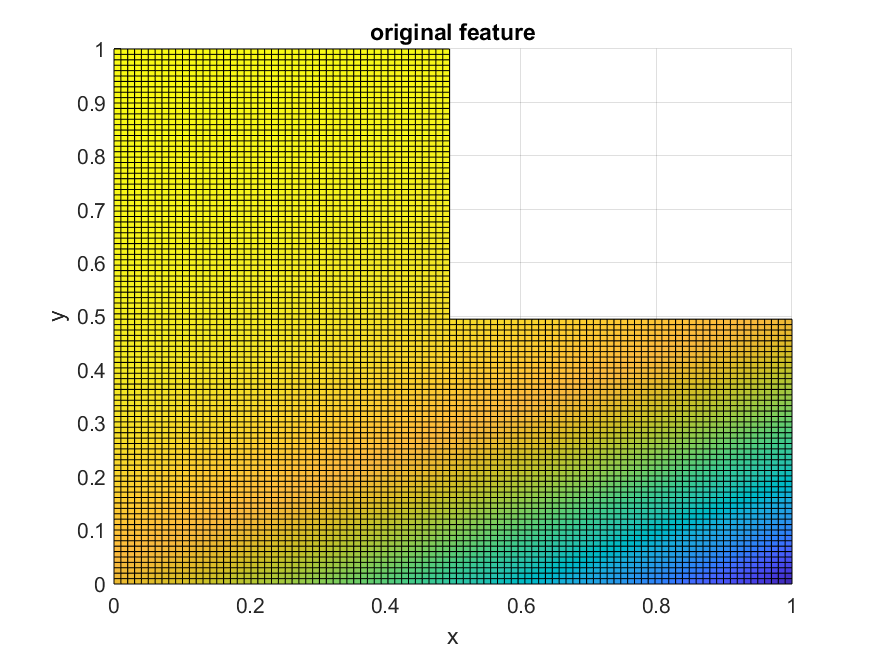}\\
        \includegraphics[width=0.45\textwidth]{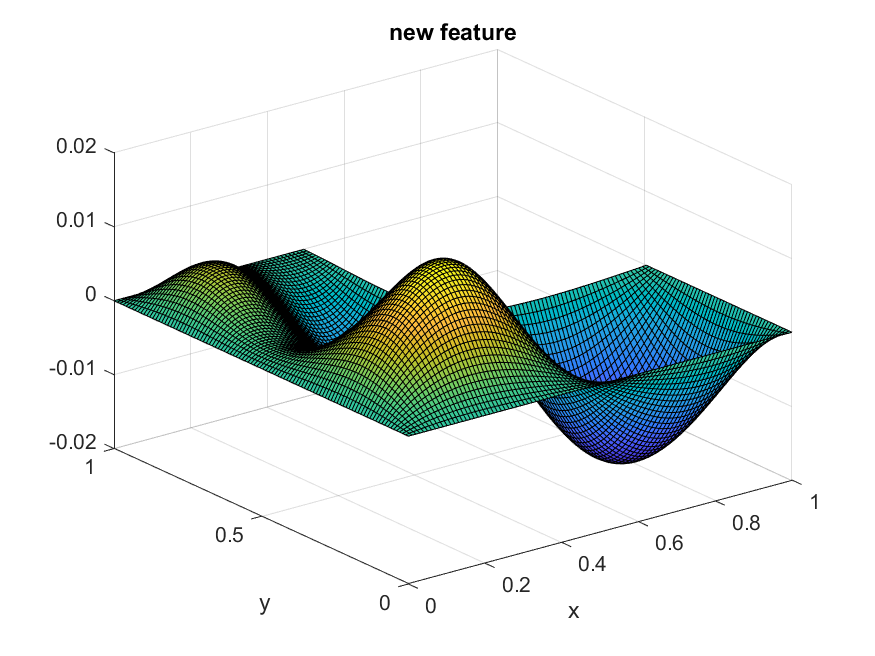}
        \includegraphics[width=0.45\textwidth]{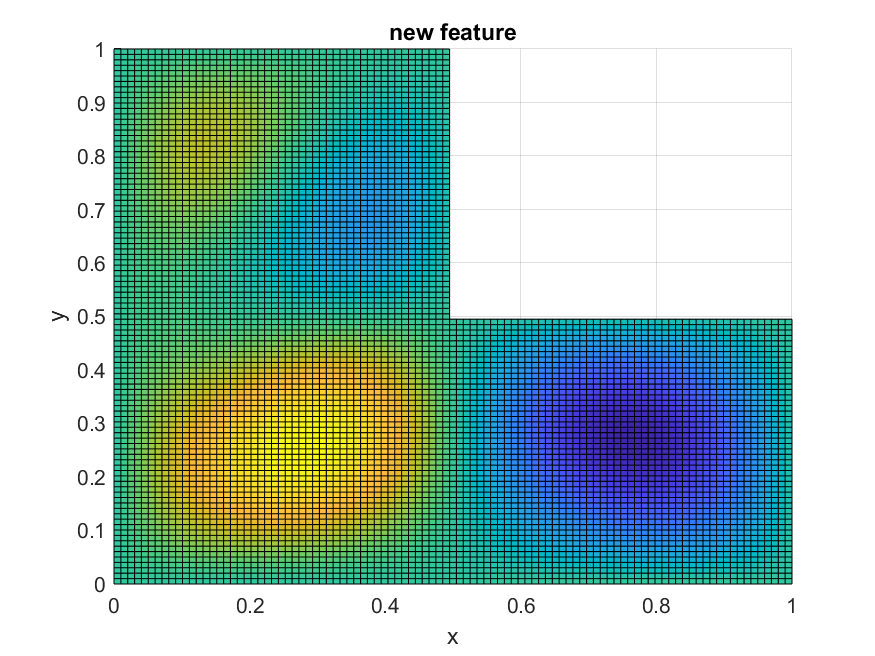}
        \caption{Top row: 3D and 2D view of the original output of a single neuron of the hidden layer (with hyperbolic tangent as activation function). It does not satisfy the boundary conditions.\\
        Bottom row: 3D and 2D view of a "projected" output. Homogeneous boundary conditions are imposed on these outputs in a least square sense with respect to the collocation points on the boundary}
    \end{figure}

    The key idea is that $\bm \beta$ is a vector in $\R^L$ and it can be uniquely decomposed as the sum of a vector in $ker(\bm{C})$ and  a vector in its complement $ker(\bm{C})^\perp$: 
    \begin{equation}
        \mathbb R^L = ker(\bm C)\, \oplus \, ker(\bm C)^\perp
    \end{equation}
    To separate the two components we can use a suitable basis, the simplest choice is a SVD decomposition:
    \begin{equation}
        \bm{C} = 
        \begin{bmatrix}
            \bm U_r \quad \bm U_\perp 
        \end{bmatrix}
        \begin{bmatrix}
            \bm \Sigma_r \quad &\bm 0 \\
            \bm 0 \quad &\bm \Sigma_\perp
        \end{bmatrix}
        \begin{bmatrix}
            \bm V_r \quad \bm V_\perp 
        \end{bmatrix}^T
    \end{equation}
    Specifically, $r$ is the numerical rank of $\bm C$, the columns of $\bm{V}_r$ are an orthonormal basis for $ker(\bm{C})^\perp$ and the columns of $\bm V_\perp$ are an orthonormal basis for $ker(\bm{C})$:
    \begin{align*}
        \bm{\beta} &= \bm{VV}^T\bm{\beta} \\
        &= \bm{V}
        \begin{bmatrix}
            \bm{V}_r^T \\ \bm{V}_\perp^T
        \end{bmatrix} \bm{\beta} \\
        &= \bm{V}
        \begin{bmatrix}
            \bm{y} \\ \bm{z}
        \end{bmatrix}
    \end{align*}
    This can be interpreted as a change of coordinates that separates the two components of $\bm \beta$. Now $y$
    can be determined using only the equality constraints:
    \begin{align*}
        \bm{C\beta} &= \bm g \\
        \bm{U \Sigma V}^T \bm{\beta} &= \bm g \\
        \begin{bmatrix}
            \bm{U}_r \quad \bm{U}_\perp
        \end{bmatrix} 
        \begin{bmatrix}
            \bm{\Sigma}_r  &\bm{0} \\
            \bm{0} &\bm{0}
        \end{bmatrix}
        \begin{bmatrix}
            \bm{y} \\ \bm{z}
        \end{bmatrix} \bm{\beta} &= \bm g \\
        \bm U_r \bm \Sigma_r\bm y &= \bm g \\
        \bm{y} &= \bm{\Sigma}_r^{-1}\bm U_r^T\bm g
    \end{align*}
    Meanwhile $\bm z$, irrelevant with respect to the equality constraints, can be determined by projecting everything on $ker(\bm C)^\perp$:
    \begin{align*}
        \bm{A\beta} - \bm{f} &= \bm{AVV}^T\bm{\beta} - \bm{f} \\
        &= \bm{A}
        \begin{bmatrix}
            \bm{V}_r & \bm{V}_\perp
        \end{bmatrix}
        \begin{bmatrix}
            \bm{y} \\ \bm{z}
        \end{bmatrix} - \bm{f} \\
        &= \bm{AV}_\perp\bm{z} - \left( \bm{f} - \bm{V}_r\bm{y} \right) \\
        &= \bm{\Tilde{A}z} - \bm{\Tilde{f}}
    \end{align*}
    A visual interpretation of this projection is given in figure \ref{fig:projected_outputs}.
    Finally, we minimize the residual in the interior of the domain:
    \begin{align}
        \bm z^* &= \arg\min \frac{1}{2}\|\bm{\Tilde{A}z} - \bm{\Tilde{f}}\|_2^2 \\
        \bm \beta^* &= \bm V
        \begin{bmatrix}
            \bm y \\ \bm z^*
        \end{bmatrix}
    \end{align} 
  This approach can be extended to nonlinear problems as follows:
    \begin{enumerate}
        \item Consider any nonlinear optimization algorithm that uses the Jacobian of the loss function
        \item At each step, compute the Jacobian with respect to $\bm z$ by the chain rule:
            \begin{equation}
                \Tilde{\bm J} = \frac{\partial \bm{\mathcal R}}{\partial \bm \beta} \frac{\partial \bm \beta}{\partial \bm z} = \frac{\partial \bm{\mathcal R}}{\partial \bm \beta} \bm V_\perp
            \end{equation}
        \item Compute the update $\bm{\delta z}$ according to the chosen algorithm
        \item Compute the updated set of parameters from  $\bm{\beta + \delta \beta} $ where:
            \begin{equation}
                \bm{\delta \beta} = \bm V_\perp \bm{\delta z}
            \end{equation}
    \end{enumerate}
    \clearpage

\section{Comparisons}\label{sect:COMPARISON}
\begin{figure}
        \centering
        \includegraphics[width=0.5\linewidth]{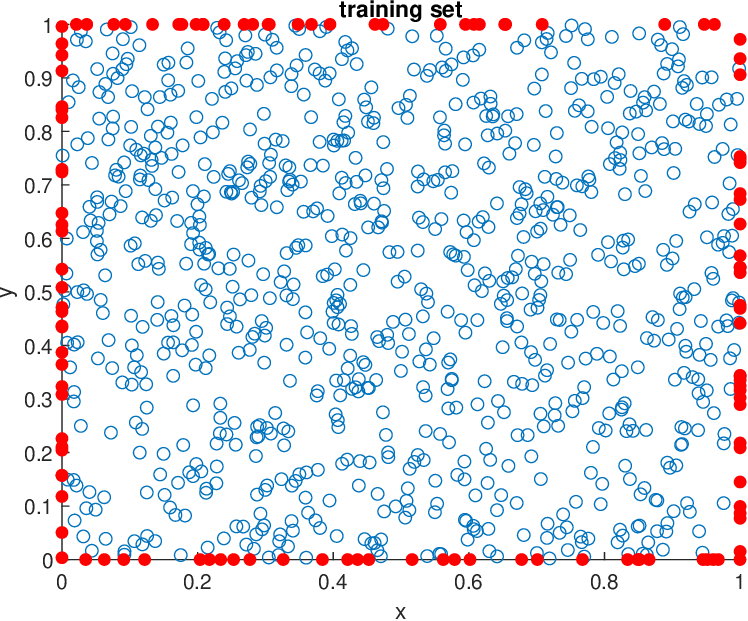}
        \caption{The domain $\Omega$ with the location of collocation nodes, also denoted by training set: internal points $\{\bm{x}_i\}_{{1 \leq i\leq N_I}}$ in blue circles and the boundary ones $\{\bm{x}_k\}_{{1 \leq k\leq N_B}}$ in red dots. In this plot we have $N=N_I+N_B=1000$. }
        \label{fig:square_training_set}
    \end{figure}

\begin{figure}
        \centering
        \includegraphics[width=0.45\linewidth]{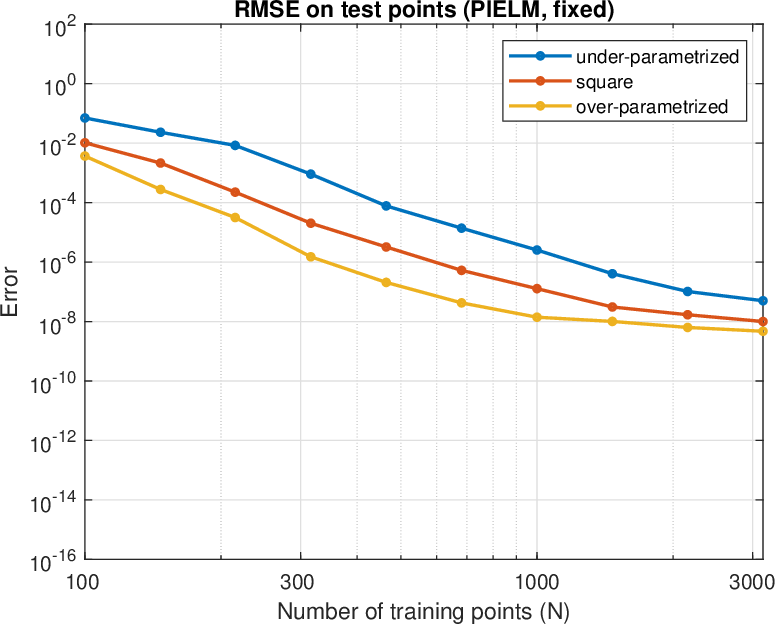}
        \includegraphics[width=0.45\linewidth]{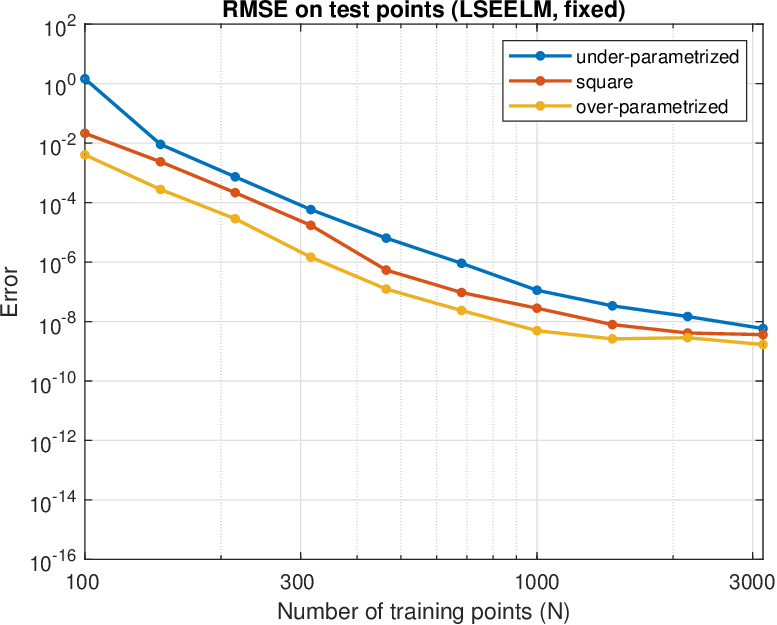}  \\
        \includegraphics[width=0.45\linewidth]{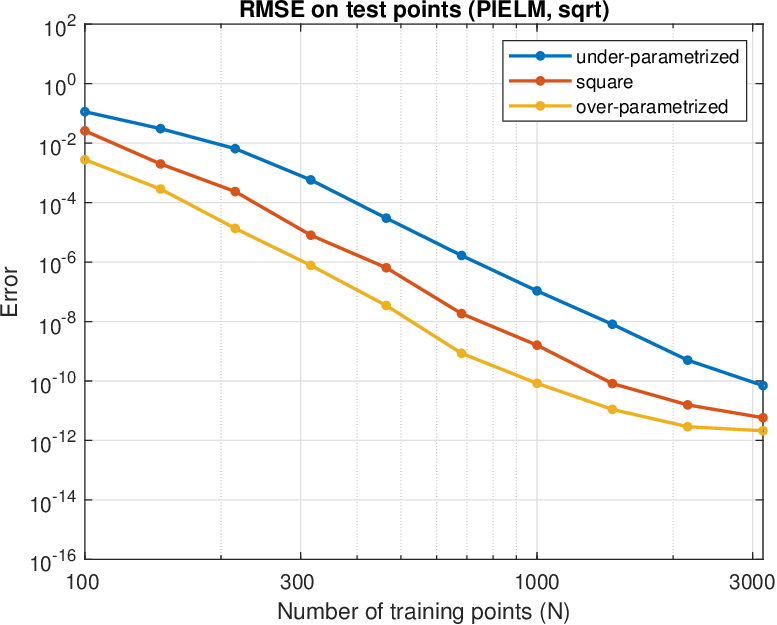}
        \includegraphics[width=0.45\linewidth]{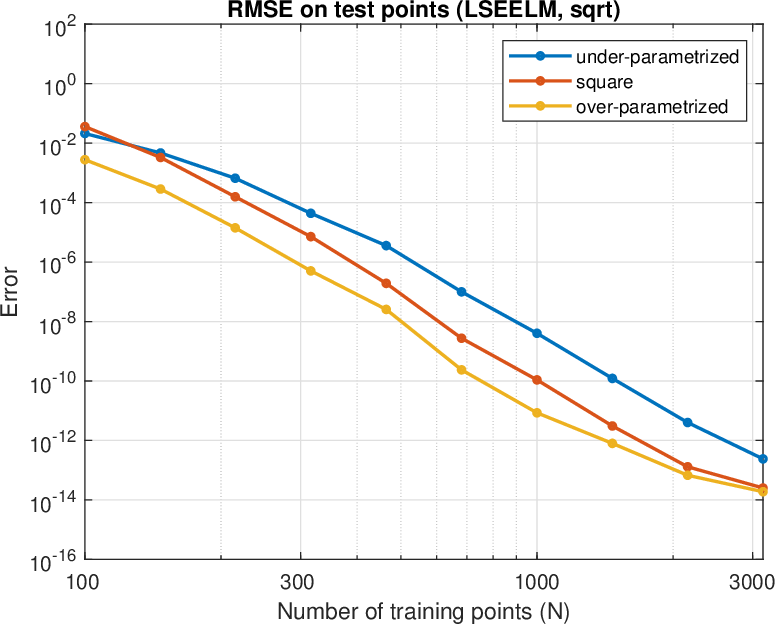}  \\
        \includegraphics[width=0.45\linewidth]{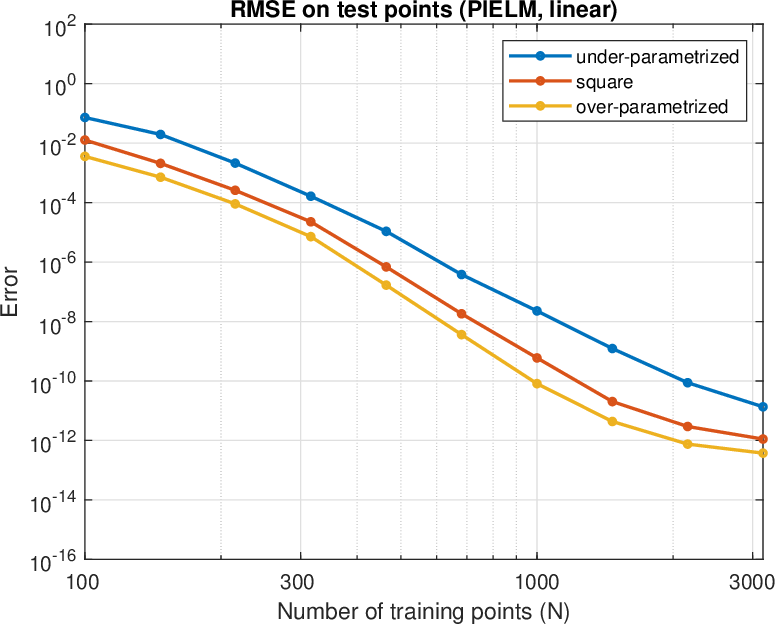}
        \includegraphics[width=0.45\linewidth]{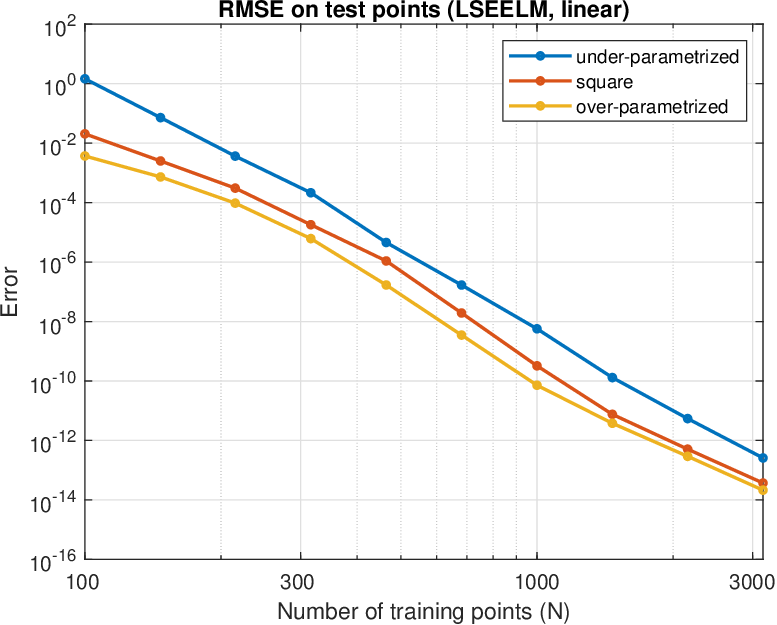}
                \caption{ Convergence results for PIELM (on the left) and LSE-ELM (on the right).\\
                Top row: with 50 training points on the boundary, and increasing only the number of points in $\Omega$. In this case we reach a plateau, with the same qualitative behavior for both methods.\\
                Central row: following the distribution $\frac{N_B}{4l}=\sqrt{\frac{N_I}{l^2}}$. The performance gets better in all cases, but we can see that PIELM underperforms in the high precision regime.\\
                Bottom row: using more points on the boundary by following the distribution $N_B = N_I = N/2$. This strategy does not lead to any improvement w.r.t the sqrt case, while using more points.}
        \label{Fig:Conv_1}
    \end{figure}

\begin{figure}
        \centering
        \includegraphics[width=0.45\linewidth]{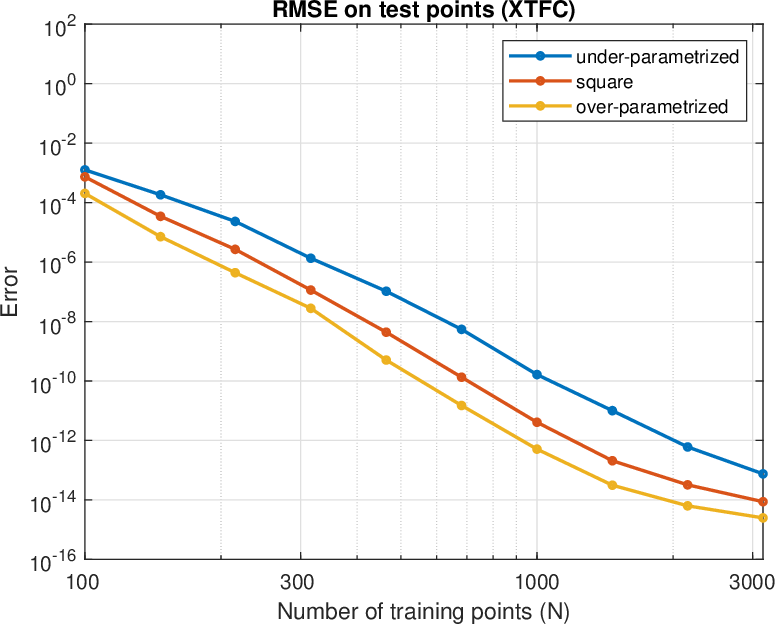}
        \caption{Convergence results for XTFC. XTFC satisfies the boundary conditions by construction, thus we only consider the scenario where $N_I=N$. Similarly to the other methods, more neurons always lead to better performances.}
       \label{Fig:Conv_2}
    \end{figure}

Having introduced in the previous section the use of ELMs in different settings, in this section tha aim is to compare such methods, LSE-ELM, PIELM, and XTFC on the 2D Poisson problem:
    \begin{equation}\label{eq:pb_sec4}
        \begin{cases}
            \Delta u(\bm x) = f(\bm x) \qquad &\forall \bm x \in \Omega \\
            u(\bm x) = g(\bm x) \qquad &\forall \bm x \in \Gamma
        \end{cases}
    \end{equation}
    where $\Omega$ is the unit square $[0,1] \times [0,1]$ and $\Gamma$ is the boundary $\partial \Omega$.
    As shown in Figure \ref{fig:square_training_set}, we take collocation points $\{\bm{x}_i\}_{{1 \leq i\leq N_I}}$ in $\Omega$ and $\{\bm{x}_k\}_{{1 \leq k\leq N_B}}$ on $\Gamma$, all sampled according to the uniform distribution.\\
    We fix the problem \ref{eq:pb_sec4} so to have a solution with high regularity. In particular, we take, see Figure \ref{fig:square_plots}: 
    \begin{equation}
        u^*(x,y) = e^{-10((x-x_c)^2+(y-y_c)^2)} + \frac{x}{10} + \frac{y}{5} \ ,
    \end{equation}
    with $x_c=y_c=0.4$. To evaluate accuracy, we take the behavior on 10000 test points, sampled randomly as done for the training set, but different from those. The reported error will be the Root Mean Squared Error (RMSE) :
    \begin{figure}
        \centering
        \includegraphics[width=0.48\linewidth]{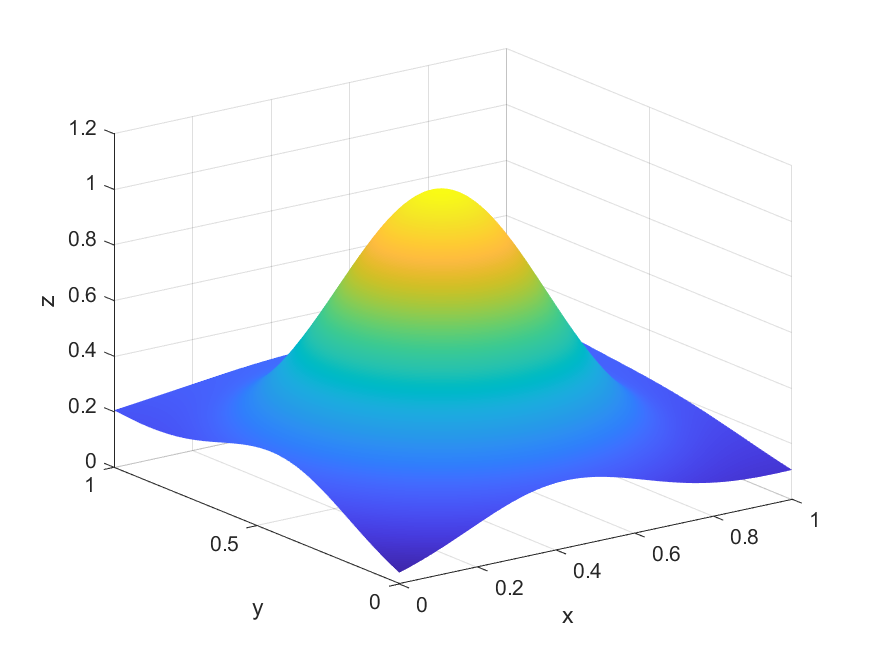}
        \includegraphics[width=0.48\linewidth]{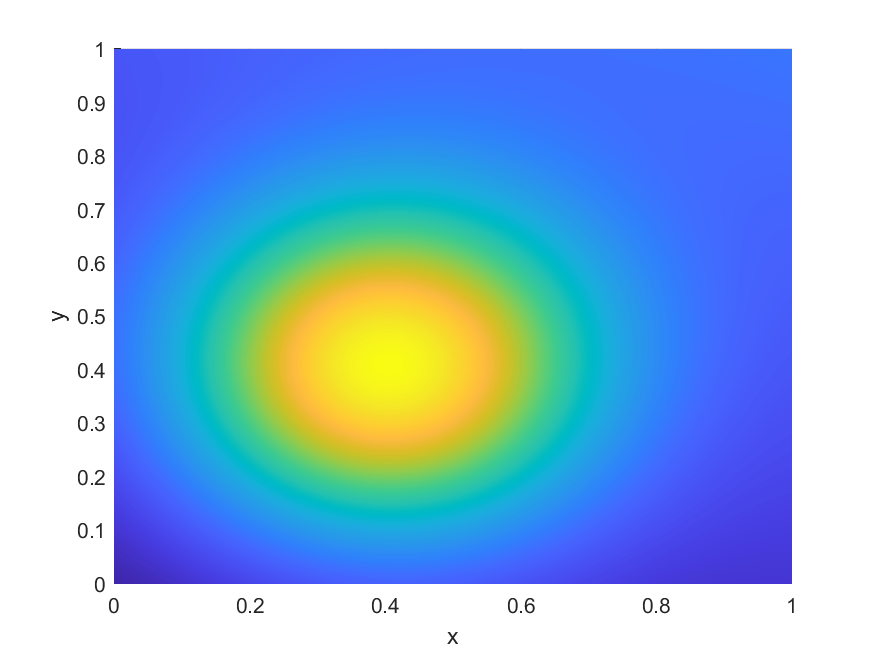}
        \caption{Exact solution for the problem considered in this section.}
        \label{fig:square_plots}
    \end{figure}
    
    \begin{equation}
        RMSE = \left( \frac{1}{N_{I,test}}\sum_{i=1}^{N_{I,test}} \Bigl(\Tilde{u}(\bm x_i) - u^*(\bm x_i)\Bigr)^2 + \frac{1}{N_{B,test}}\sum_{k=1}^{N_{B,test}} \Bigl(\Tilde{u}(\bm x_k) - u^*(\bm x_k)\Bigr)^2\right)^\frac{1}{2}
    \end{equation}
    We are interested in how the performance changes as some parameters of the problem are modified, the following are all the options we took into account:
    \begin{itemize}
        \item We compare 3 different ELM training strategies (LSE-ELM, PIELM, and XTFC).

        \item For PIELM and LSE-ELM, when the total number of training points $N:=N_B+N_I$ increases, we considered 3 options for the proportion between $N_I$ and $N_B$:
            \begin{enumerate}
                \item $N_B$ is fixed and only $N_I$ increases, referred as "fixed".
                \item $N_B = N_I = N/2$, referred as "linear".
                \item The linear density of training points on $\Gamma$ is the square root of the surface density of training points in $\Omega$: $\frac{N_B}{4l}=\sqrt{\frac{N_I}{l^2}}$, referred as "sqrt".
            \end{enumerate}
                    \item We considered 3 types of convergences:
            \begin{enumerate}
                \item $N$ is fixed and $L$ increases.
                \item $L$ is fixed and $N$ increases.
                \item Both $N$ and $L$ increase keeping a constant ratio.
            \end{enumerate}
        \item When both $N$ and $L$ increase, we considered 3 ratios:
            \begin{enumerate}
                \item An under-parametrized network, $L=N/2$, referred as "under-parametrized".
                \item A network with as many neurons as training points, $L=N$, referred as "square".
                \item An over-parametrized network, $L=2N$, referred as "over-parametrized".
            \end{enumerate}
    \end{itemize}
    All convergence results are presented in Figures \ref{Fig:Conv_1}-\ref{Fig:Conv_2}. In the relative captions, we discuss the behavior. We remark that all the procedures are trained in randomly generated points and with randomly generated ELM networks, and that all could benefit from some adaptive procedure. However, in order to present general strategies, in the present paper, we focus on such setting and will consider such adaptivity strategies in future work. 

\begin{figure}
        \centering
        \includegraphics[width=0.45\linewidth]{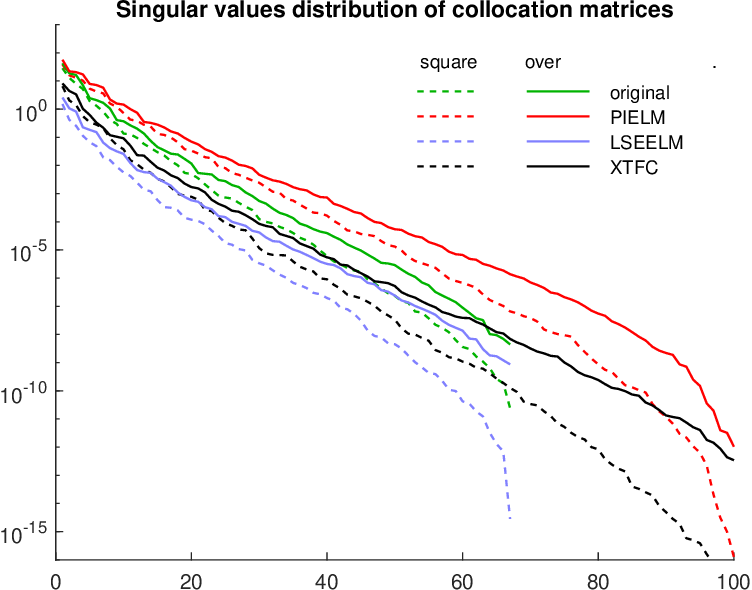}
        \includegraphics[width=0.45\linewidth]{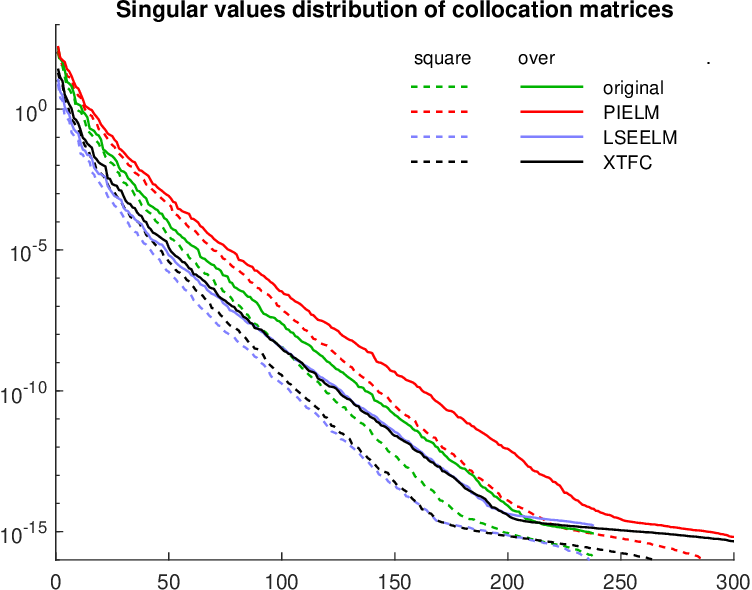} \\
        \includegraphics[width=0.45\linewidth]{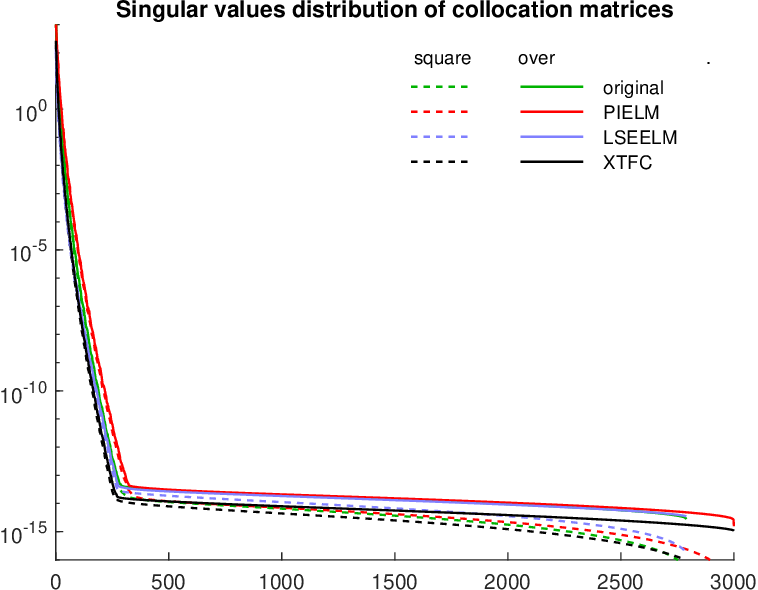}
        \includegraphics[width=0.45\linewidth]{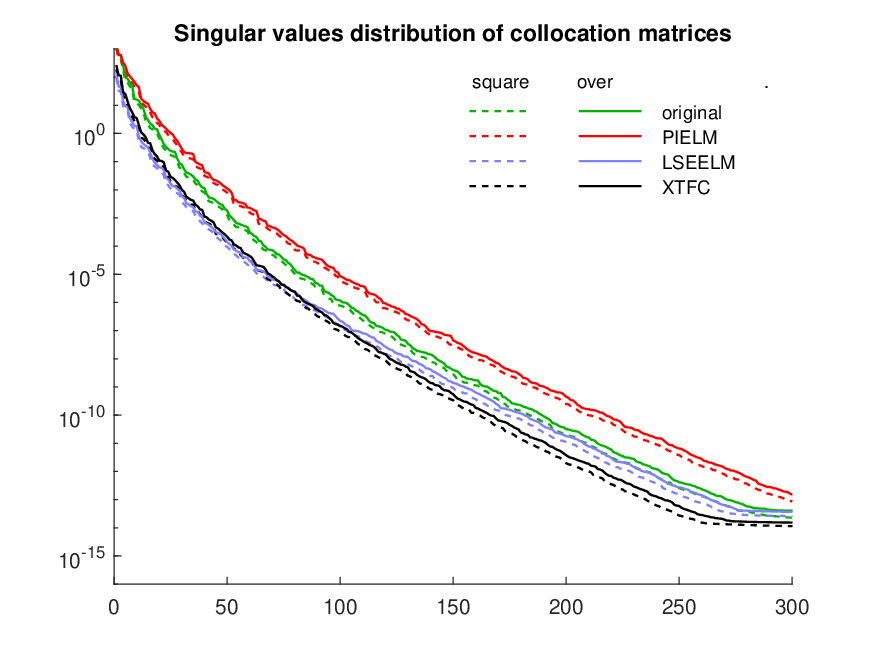}
        \caption{Singular values distributions for the same problem, but with training sets and network of different sizes. "Square" refers to a network with as many neurons as training points, while "Over" refers to a network with twice as many neurons as training points. \\
        Top row: 100 (left) and 300 (right) training points, the "original" collocation matrix and the one computed in "LSE-ELM" only use internal points, therefore their spectrum is shorter. \\
        Bottom row: 3000 training points (left), cut to the first 300 (right). While the "numerical" rank of the matrix does not grow with the size of the training set, we remark that from the previous convergence results however the resulting approximation still gets better.}
        \label{fig:singular_tutti}
    \end{figure}

    As a general comment, we can notice that the accuracy for the three procedures is similar, but some remarks on the computational costs can be highlighted. PIELM solves a unique linear system of dimension $N\times L$, XTFC also solves a unique system of the same type, but we remark that in this case $N=N_I$ and a pre-processing for the resolution of the boundary problem is needed. The proposed procedure, on the other hand, solves two linear systems, one of dimension $N_B\times L$ and the final one with $N_I$ rows. This gives a saving in terms of computational cost and has an impact on the condition number of the overall resolution. Following what was done in \cite{Chen2024}, we present in Figure \ref{fig:singular_tutti} the comparison between the spectra obtained via SVD. Fixing the dimension $N=100, 300, 3000$, in the plots we present both the square and the $L=2N$ overparameterized case. Different colors refer to different matrices: the collocation of the original ELM functions on the $N_I$ internal points (denoted by "original"), the collocation on all $N_I+N_B$ points for the PIELM method, the collocation of the modified functions on the $N_I$ internal points for the LSE-ELM method, the collocation on $N$ (all internal) points for the XTFC method. From the plot, it appears that the spectra reach machine precision after $\approx 250$ values, we remark that this could depend on the very simple structure both of the problem and of the domain.

The final conclusion that we can draw from these comparisons is that overparametrizzation helps in terms of accuracy and, mostly, in terms of conditioning of the problem. Our new procedure takes great advantage from the subdivision of the problem in two sub-problems so to have an overall condition number that is order of magnitude lower than the direct competitors, and in all cases is very similar to the collocation on original points that does not consider the boundary conditions. In the next numerical tests, moreover, we consider the sqrt case for the balancing of internal and boundary points.


    \graphicspath{{./Figures}}

\section{Numerical tests}\label{sect:TEST}

In this section, we present our numerical tests. First of all, we consider two simple problems with known solution in general domains and evaluate the distribution of errors, thus accuracy. In order to present the good generalization properties of our method, we also present results where the additional conditions on $\Gamma$ are internal conditions.

Then, we consider the convergence behavior of the proposed method in the case of challenging problems. In Section \ref{sec:peak}, we consider peaked solutions, in Section \ref{sec:boundary_layer} problems with a boundary layer, and in Section \ref{sec:Anisotropic_diffusion} anisotropic diffusion: for these prolems we notice that our procedure always converges, but the speed of convergence strongly depends on the presence of steep gradients. The last problem that we consider is the reentrant corner in Section \ref{sect:corner}: for this problem, we notice that the presence of lower regularity gives that our method fails to converge, as expected for spectral methods. 

Finaly, we consider some examples of nonlinear problems, following the approach for the resolution of the nonlinearity presented in section \ref{Sect:Constrained}. In particular two time-dependent nonlinear parabolic problems are used.

For all the problems, we fix the distribution of internal and boundary collocation points as $\frac{N_B}{4l}=\sqrt{\frac{N_I}{l^2}}$ and we consider over-parameterized networks: $L = 2N = 2(N_I+N_B)$. The collocation points are chosen randomly according to the uniform distribution in the domain $\Omega$ and on the boundary $\Gamma$. Finally, the random weights and biases in the activation functions are sampled independently form the uniform distribution on $[-3,3]$. 

In the case of time-dependent problems, the inial condition is teated as a fixed boundary condition, while no condition is imposed at the boundary of the domain corresponding to the final time.




\subsection{L-shape and star shape domain}
\begin{figure}[H]
        \centering
        \includegraphics[width=0.45\linewidth]{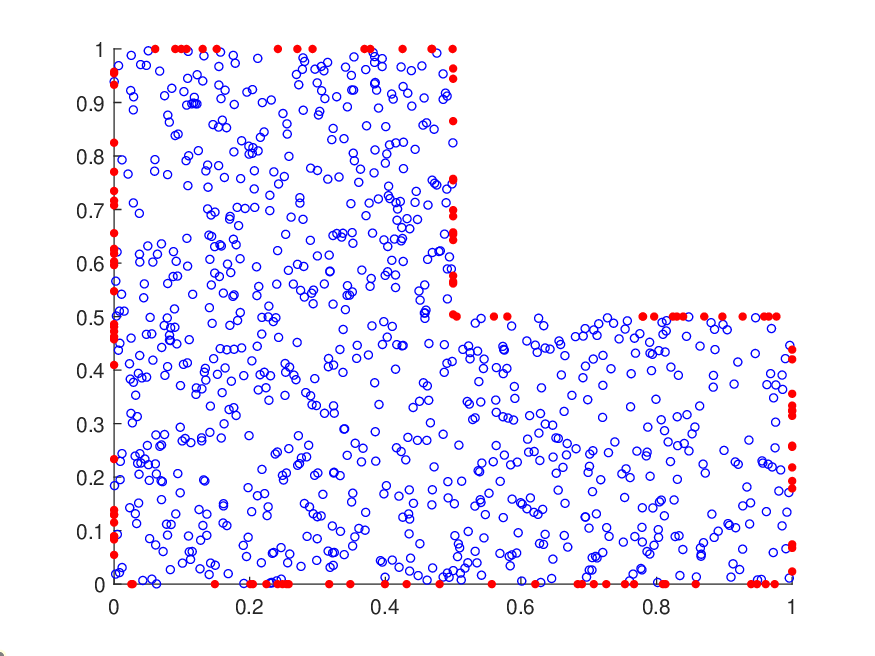}
        \includegraphics[width=0.45\linewidth]{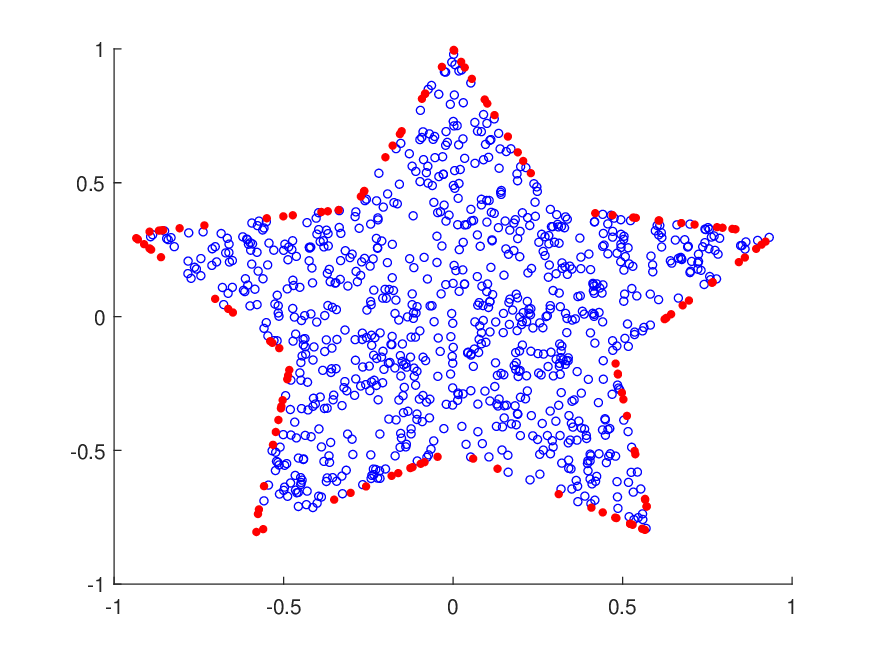}
           \caption{Training points considered, here $N_I=881, N_B=119$.}
           \label{Fig:lshape1}
    \end{figure}

\begin{figure}[H]
        \centering
        \includegraphics[width=0.45\linewidth]{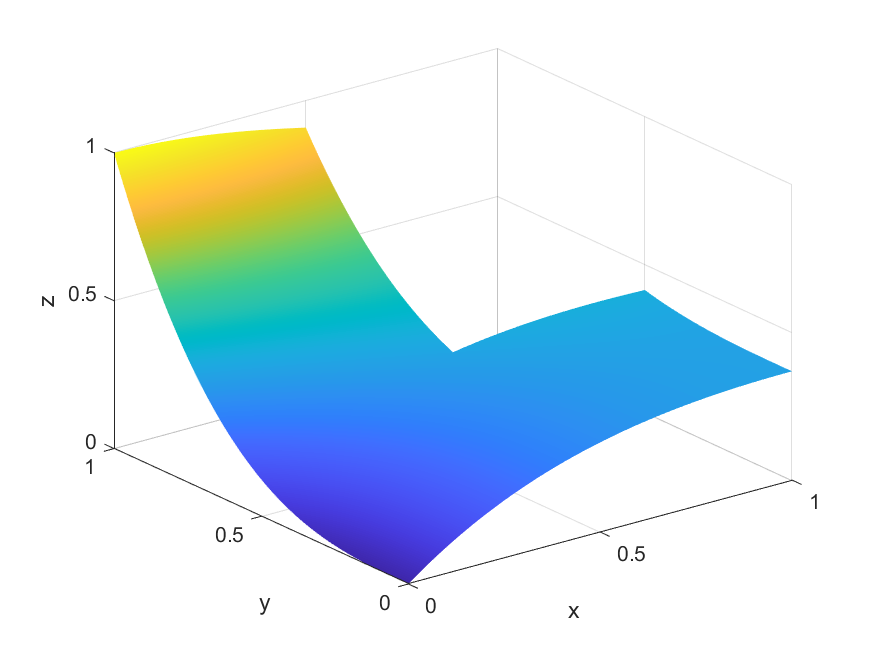}
        \includegraphics[width=0.45\linewidth]{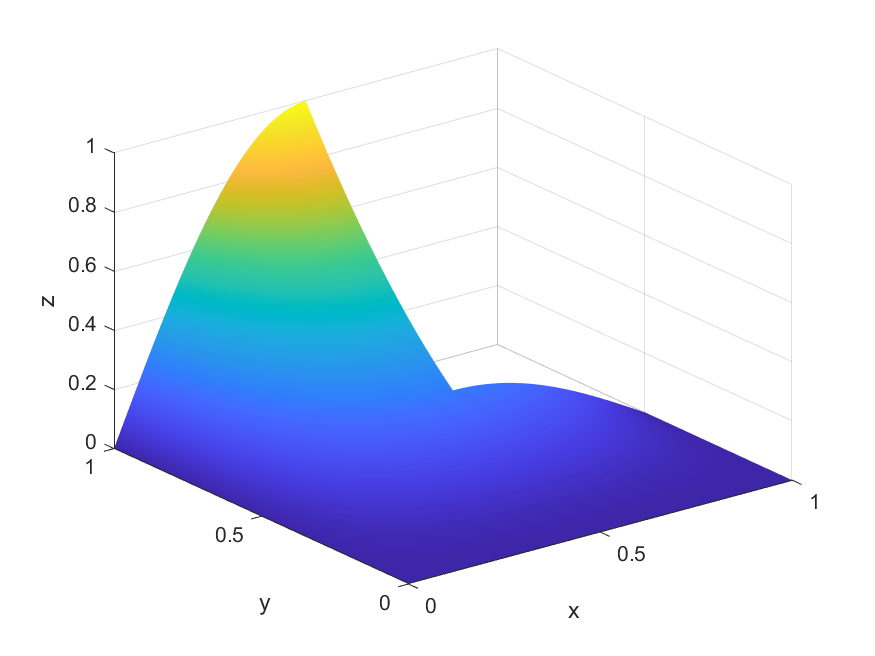}\\
         \includegraphics[width=0.45\linewidth]{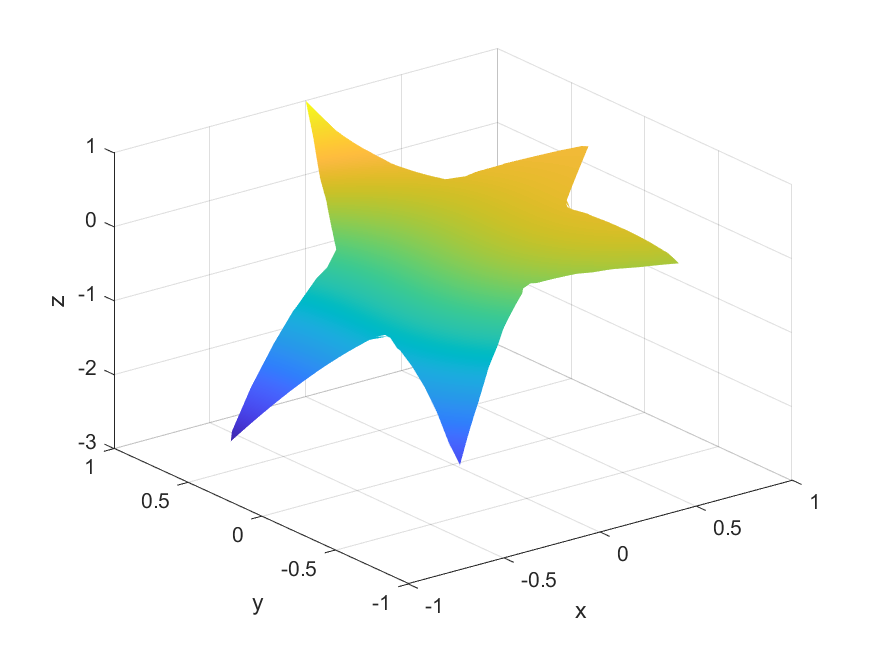}
        \includegraphics[width=0.45\linewidth]{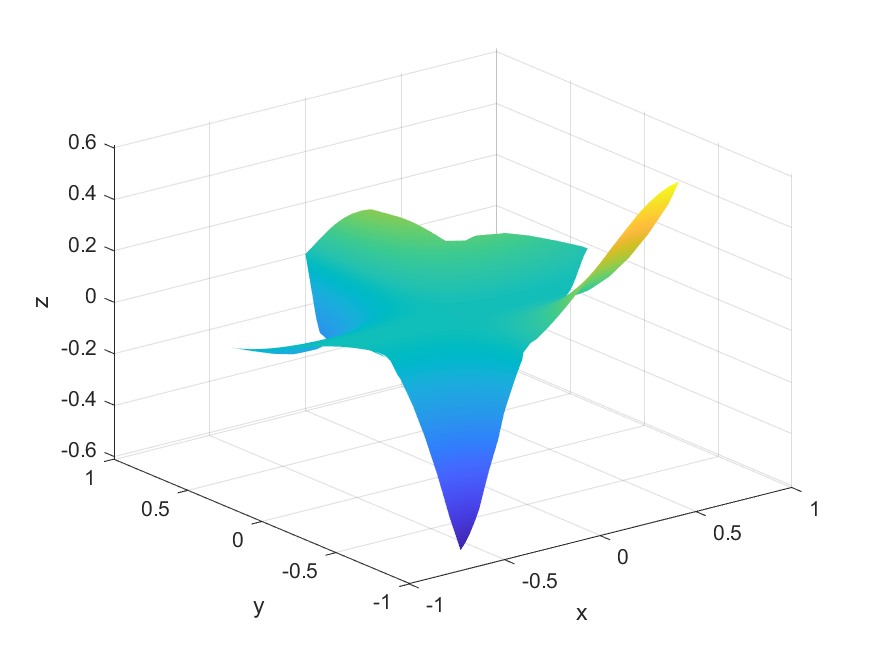}       
        \caption{Exact solutions for the proposed problems.}
\label{Fig:lshape2}        
    \end{figure}

\begin{figure}[H]
        \centering
        \includegraphics[width=0.45\linewidth]{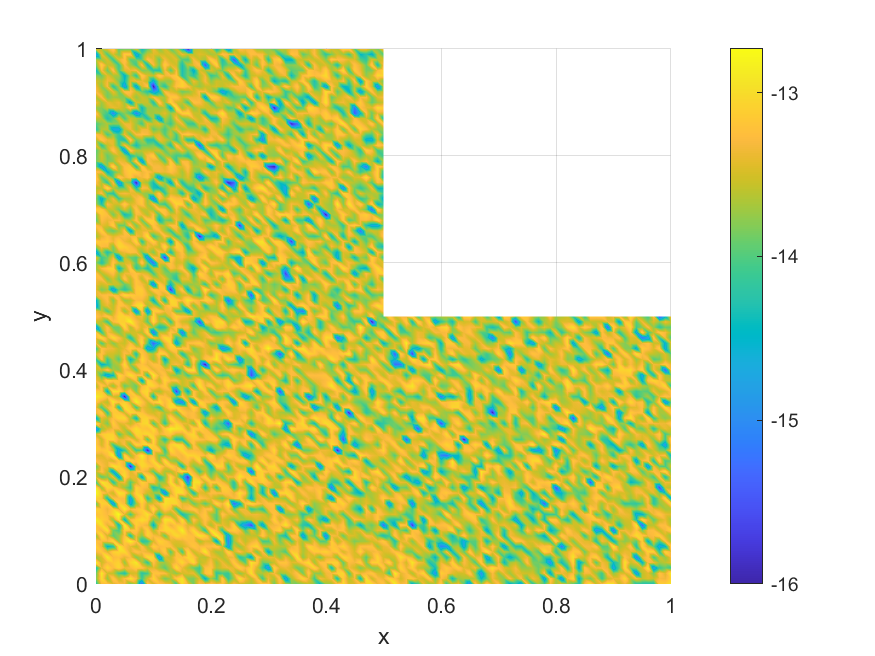}
        \includegraphics[width=0.45\linewidth]{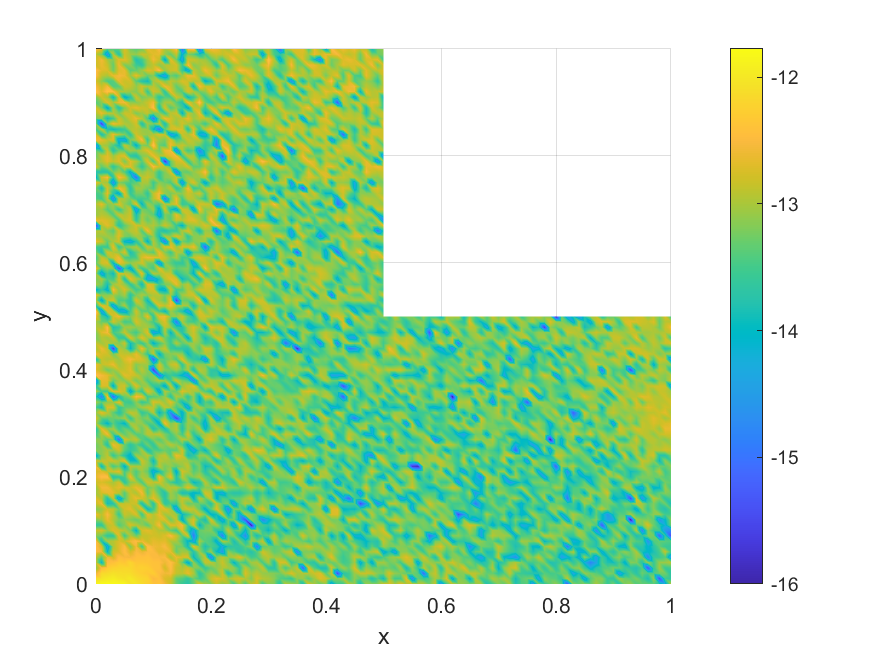} \\
        \includegraphics[width=0.45\linewidth]{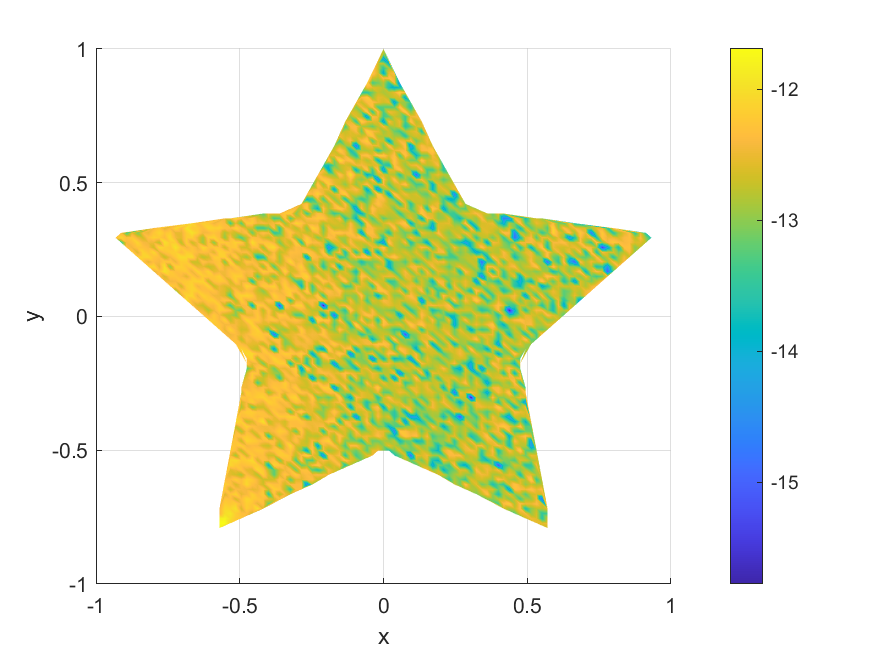}
        \includegraphics[width=0.45\linewidth]{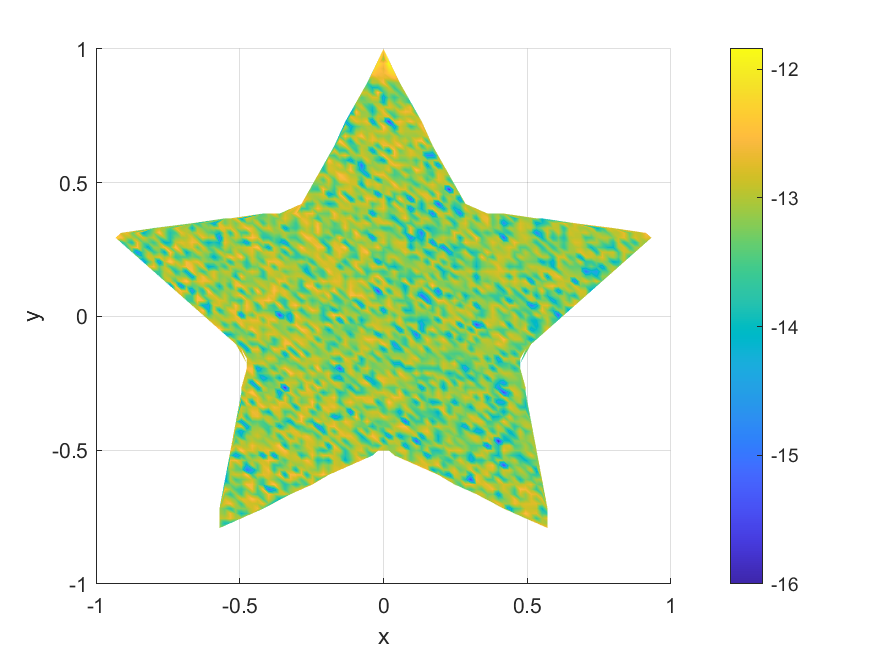} 
        \caption{Absolute errors, on the left for Problem 1 and on the right Problem 2. Notice that the values are for both problems of the order of $10^{-12}$.}
\label{Fig:lshape3}        
    \end{figure}

\begin{figure}
    \centering
    \includegraphics[width=0.48\linewidth]{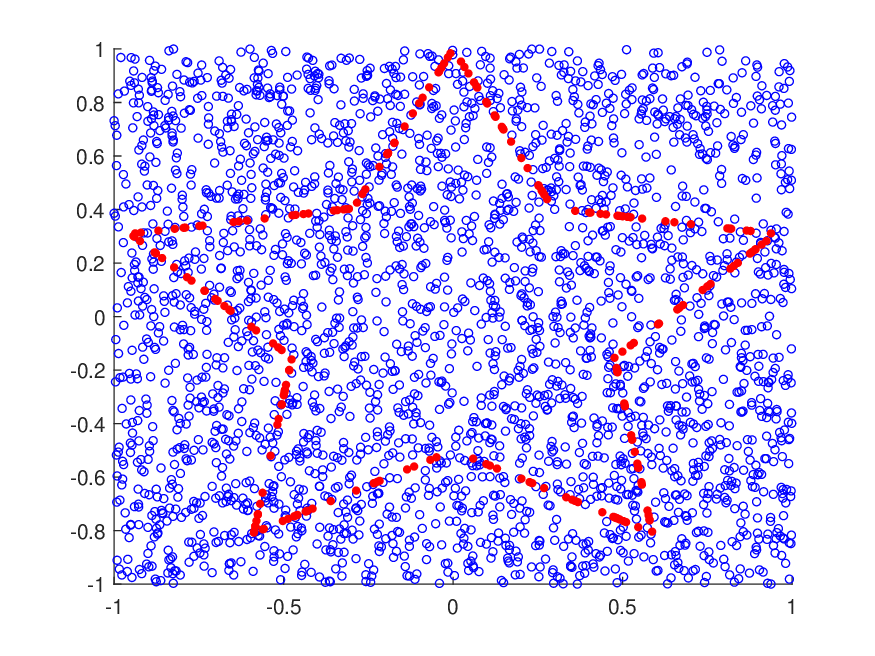}
    \includegraphics[width=0.48\linewidth]{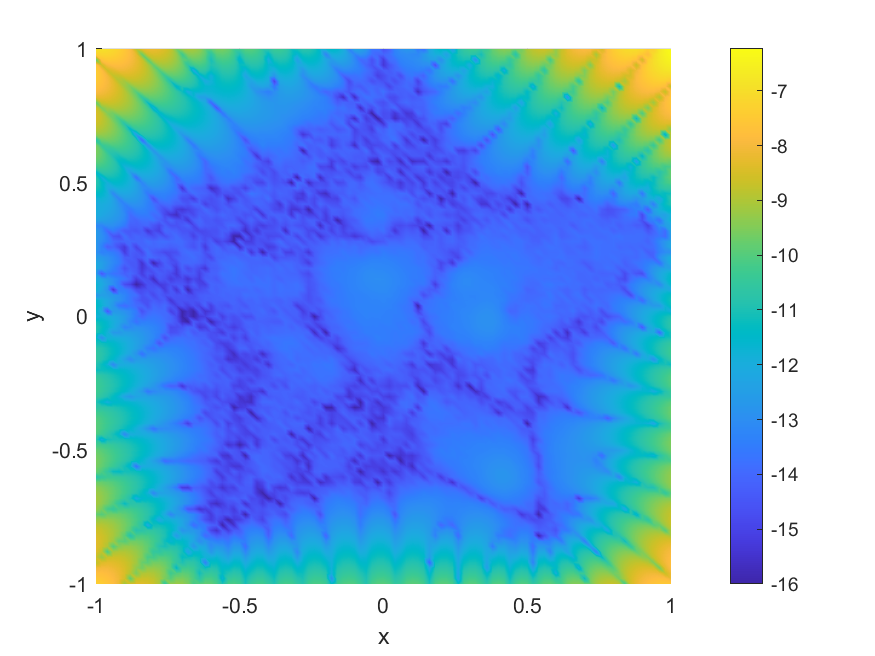}\\
    \includegraphics[width=0.48\linewidth]{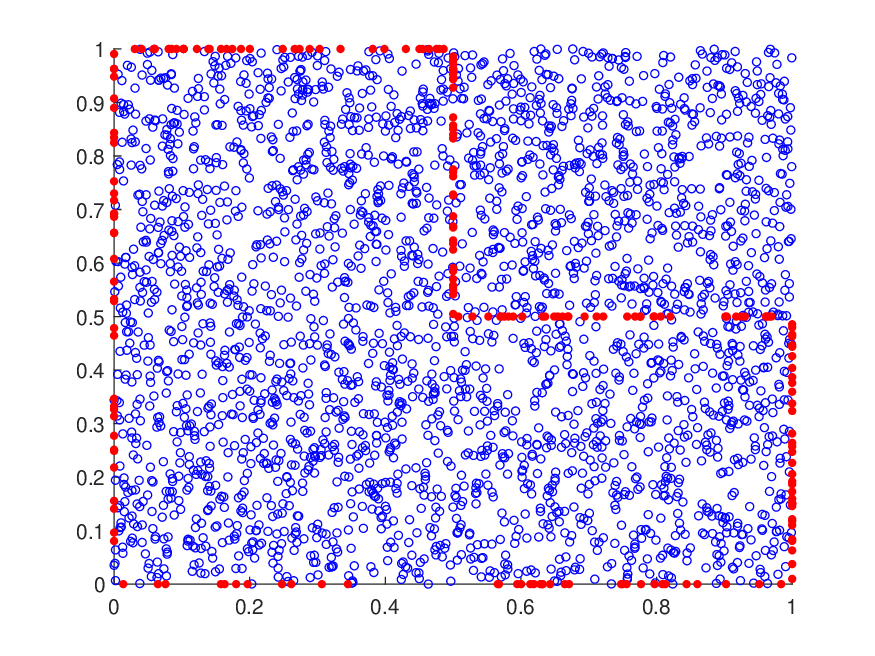}
    \includegraphics[width=0.48\linewidth]{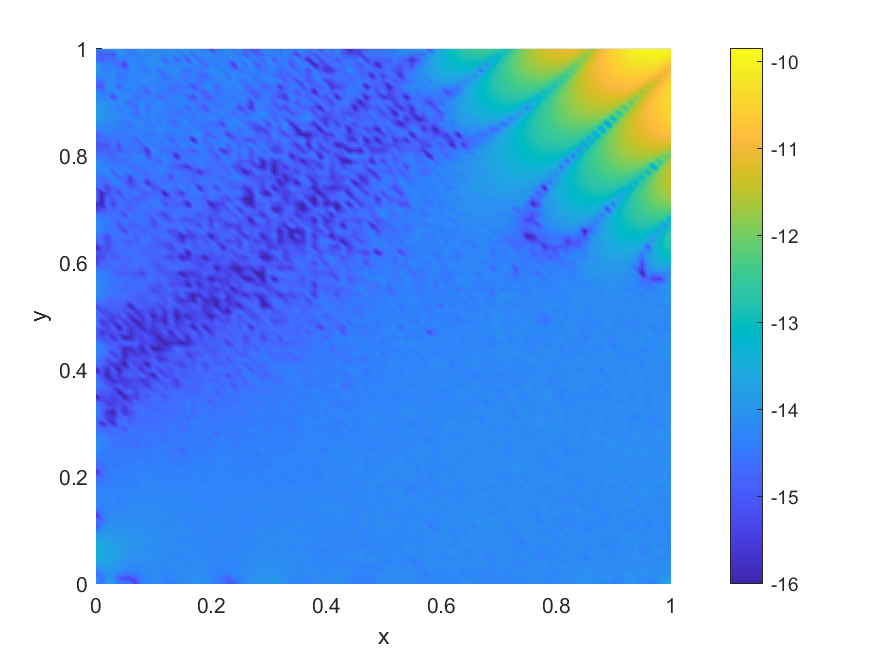}
    \caption{Absolute error in  $log_{10}$ base, extrapolation case.\\ The Poisson problem is solved in a square domain and the fixed conditions are imposed on a curve $\Gamma$ that does not coincides with the boundary. On the top, the  values are enforced on a star inside the domain. On the bottom, the values are enforced on an L-shaped region that partially intersect with the boundary of the domain. In these tests we use {$N_I = 2789, N_B = 211$}. In both cases we obtain very accurate results inside and good results also in the external part.}
    \label{fig:extrapolation}
\end{figure}

    In these first are four examples where we solve the Poisson problems in L- and star-shaped domains with Dirichlet boundary conditions to have as exact solution these two: $u(x,y)= e^{-x}(x + y^3)$ for Problem 1 and $u(x,y)= y^2 \, sin(\pi \, x)$ for Problem 2. The selected points and exact solutions are reported in Figures \ref{Fig:lshape1}-\ref{Fig:lshape2}; results are available in Figure \ref{Fig:lshape3}. At the end of this section, we consider the proposed shapes as internal curves and solve an extrapolation problem, see Figure \ref{fig:extrapolation}.


\subsection{Peak}
    \label{sec:peak}
\begin{figure}[H]
        \centering
        \includegraphics[width=0.45\linewidth]{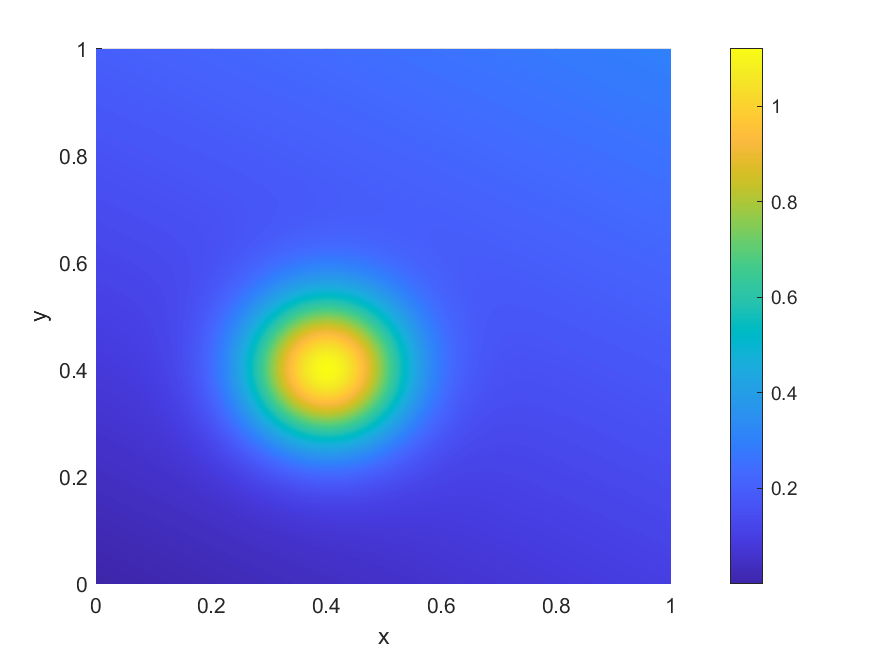}
        \includegraphics[width=0.45\linewidth]{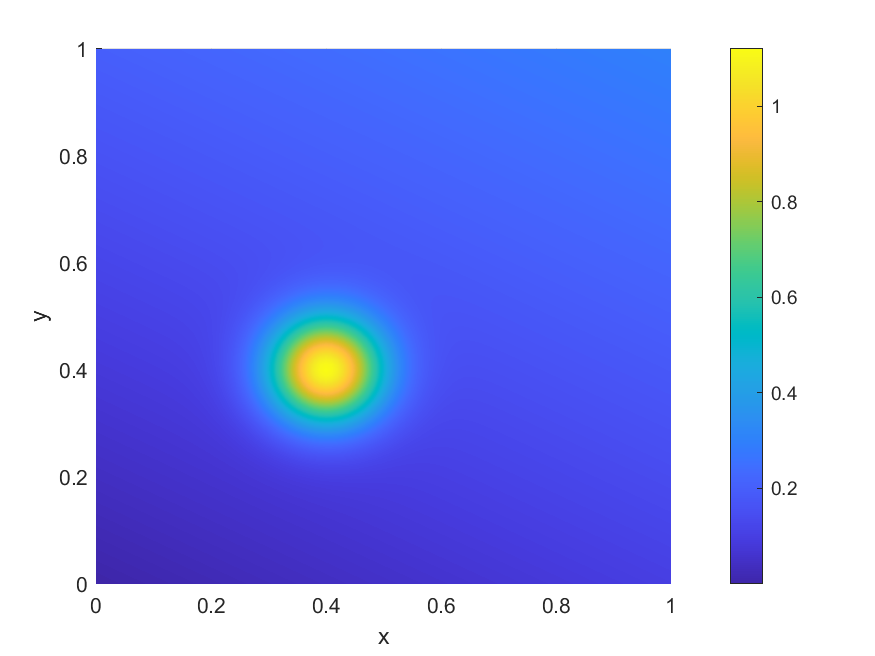}\\
        \includegraphics[width=0.45\linewidth]{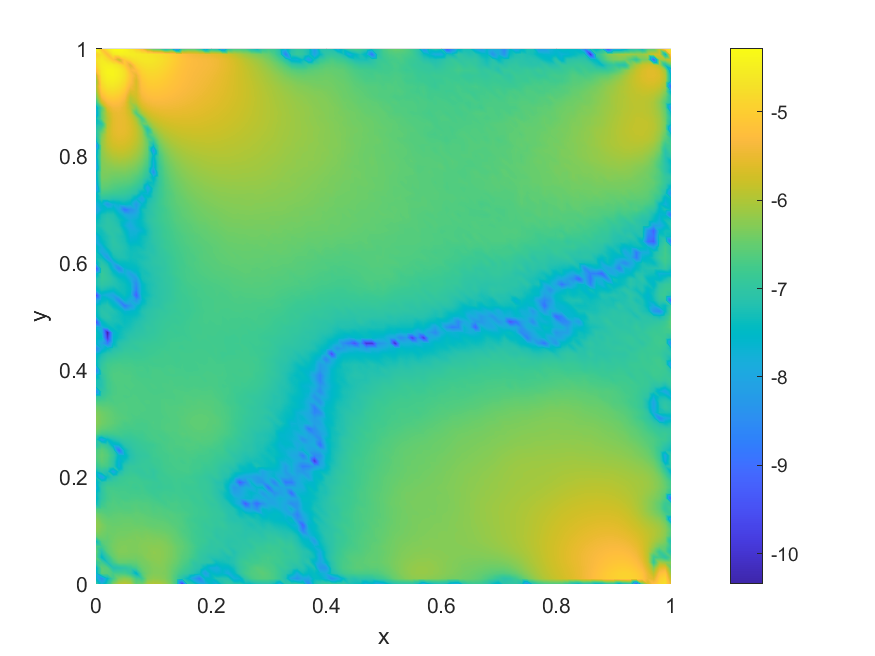}
        \includegraphics[width=0.45\linewidth]{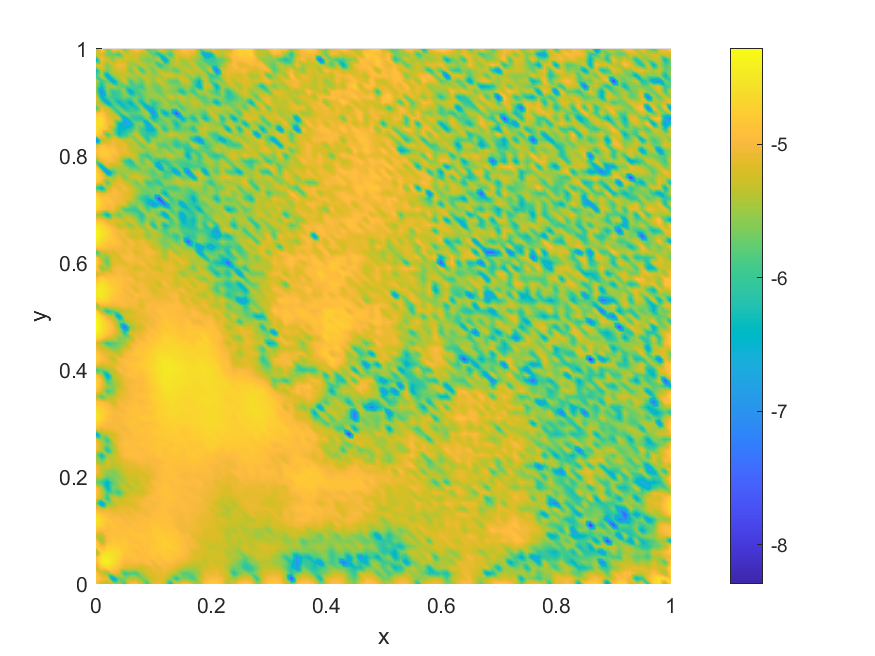}
        \caption{Poisson problem. On the left, the case where $\alpha=50$, on the right $\alpha=100$. The latter gives higher gradients. On the top, the exact solution, on the bottom computed errors in the over-parametrized case with $N=1000$.}
        \label{fig:peak}
    \end{figure}

\begin{figure}[H]
        \centering
        \includegraphics[width=0.45\linewidth]{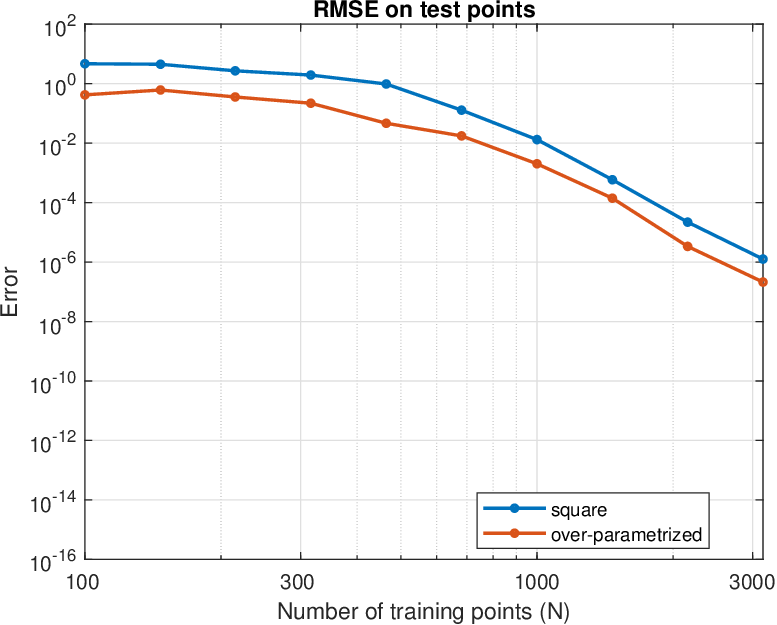}
        \includegraphics[width=0.45\linewidth]{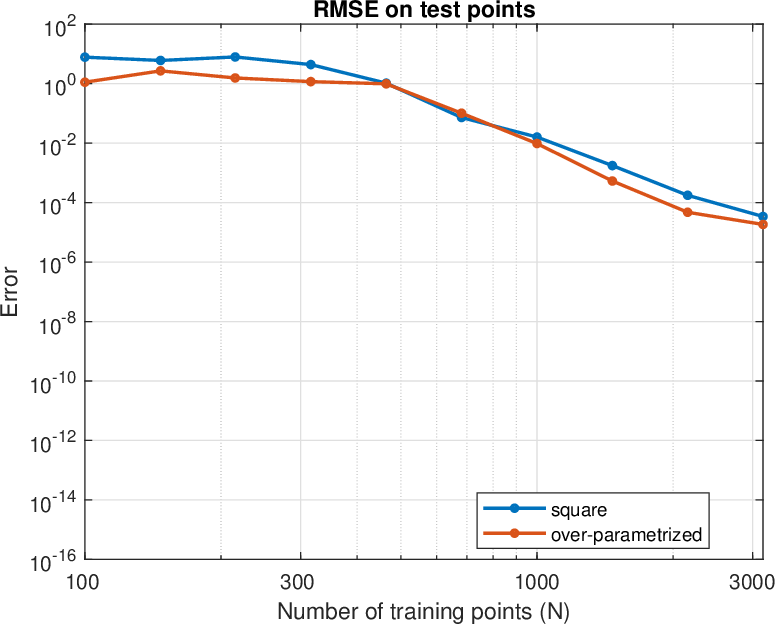}
        \caption{Convergence plots for the Poisson problems. On the left, the case where $\alpha=50$, on the right $\alpha=100$.}
        \label{fig:peak_2}
    \end{figure}

    In these examples, we consider a Poisson problem. The domain is $\Omega=[0,1]^2$ and the inside problem is $\Delta u = f$. The source term $f$ and the Dirichelet boundary conditions are fixed to have an exact solution:
    $u(x,y)=e^{-\alpha((x-x_c)^2+(y-y_c)^2)}$. We fix $x_c = y_c = 0.4$ and consider the two cases $\alpha=50, 100$, and report in Figures \ref{fig:peak}-\ref{fig:peak_2} the results obtained.


\subsection{Boundary layer}
\label{sec:boundary_layer}
\begin{figure}[H]
        \centering
        \includegraphics[width=0.45\linewidth]{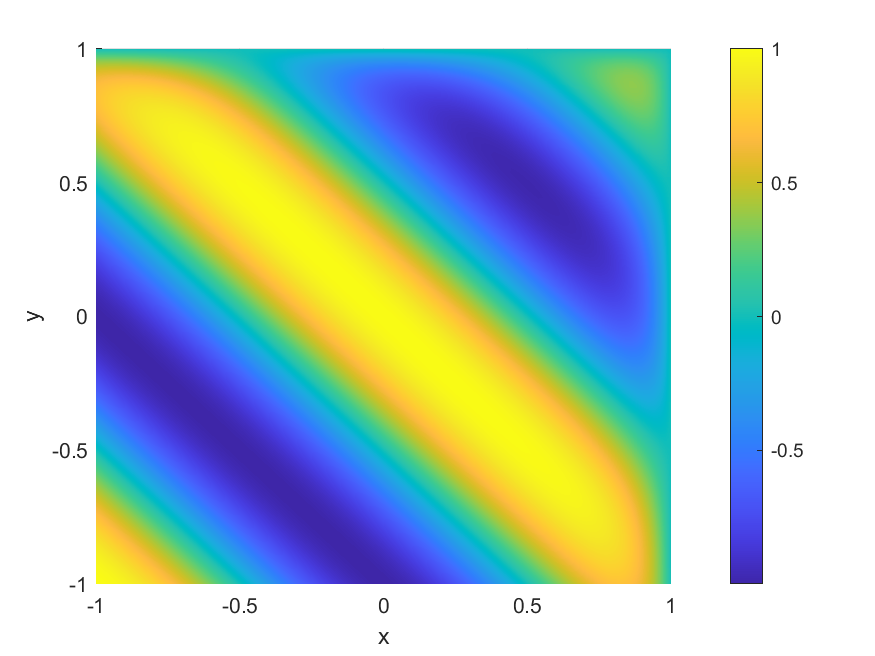}
        \includegraphics[width=0.45\linewidth]{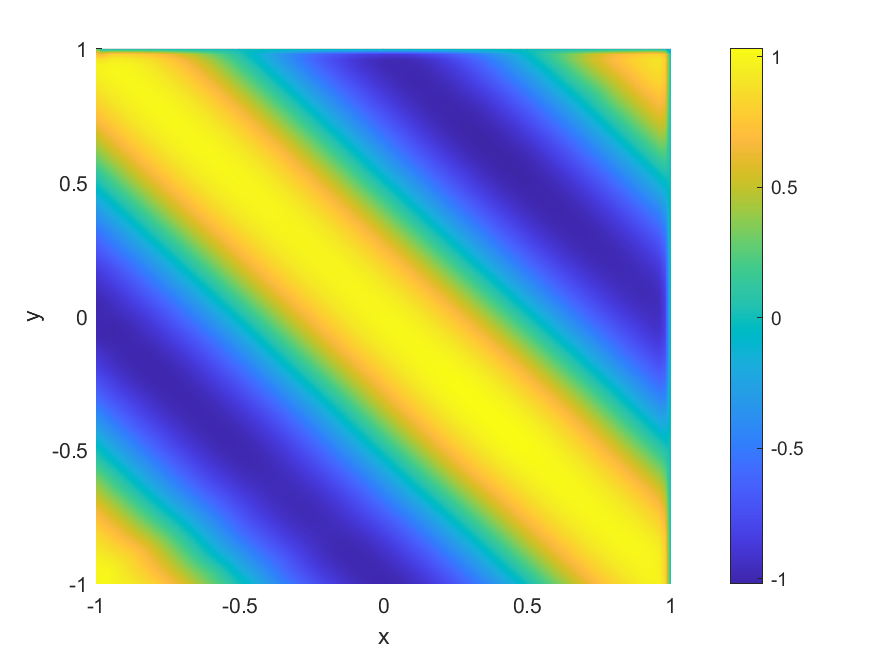}
        \\ 
        \includegraphics[width=0.45\linewidth]{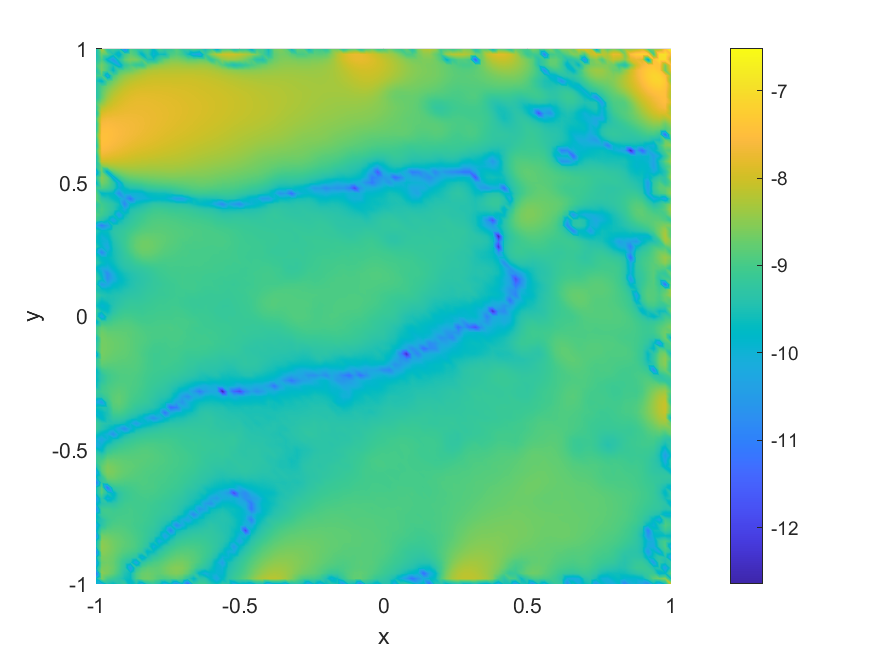}
        \includegraphics[width=0.45\linewidth]{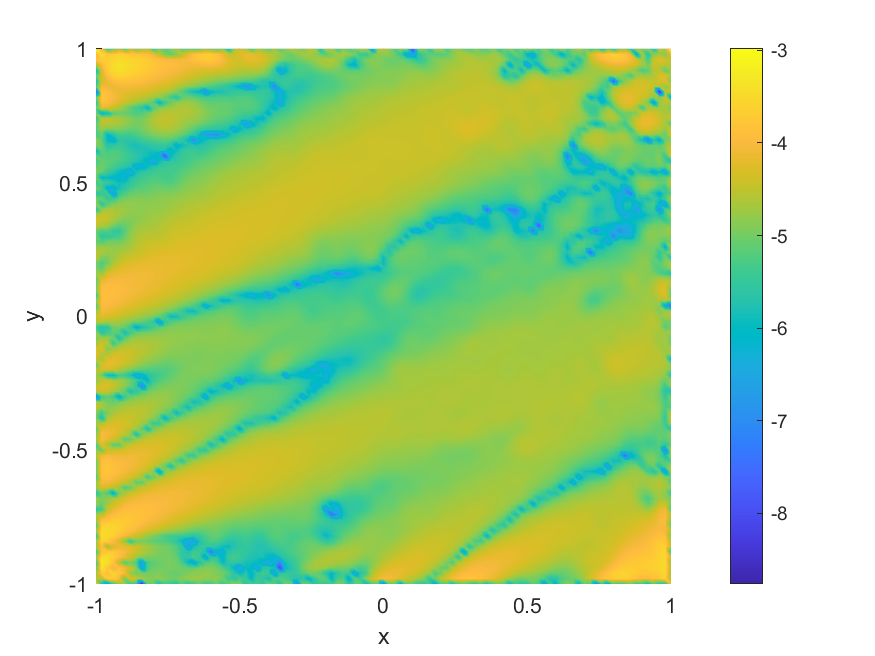}
        \caption{Boundary layer problem. On the left, the case where $\delta=0.1$, on the right $\delta=0.01$. The latter gives a stiff behavior. On the top, the exact solution, on the bottom computed errors in the case $N=1000$.}
        \label{Fig:boundarylayer}
    \end{figure}

\begin{figure}[H]
        \centering
        \includegraphics[width=0.45\linewidth]{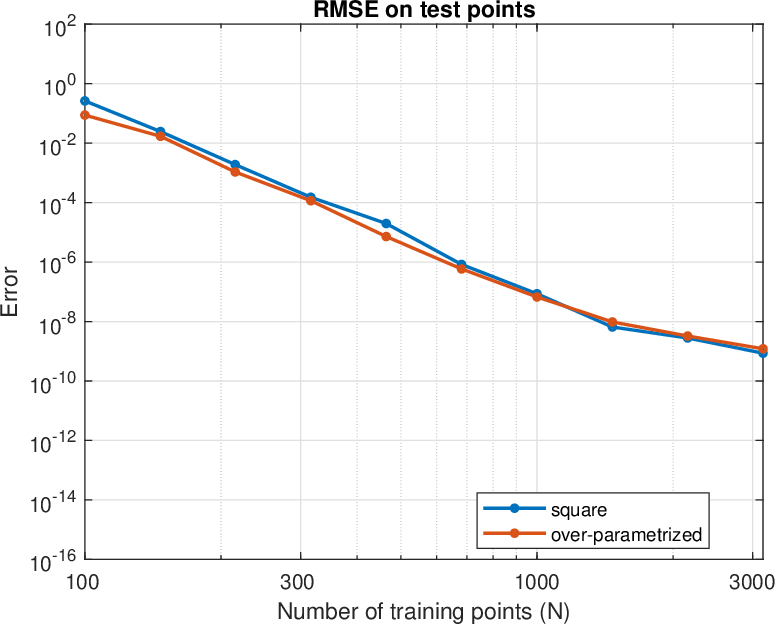}
        \includegraphics[width=0.45\linewidth]{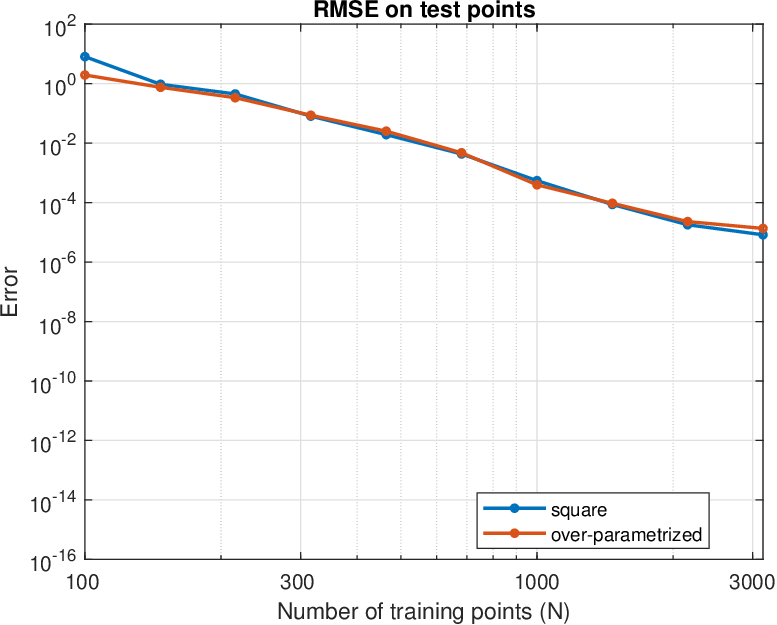}
        \caption{Convergence plots for the boundary layer problems. On the left, the case where $\delta=0.1$, on the right $\delta=0.01$. As it can be seen, the accuracy strongly depends on the stiff behavior of the solution.}
        \label{Fig:boundarylayer_2}
    \end{figure}

In these examples, we consider diffusion-trasport problems that pose a boundary layer. The domain is $\Omega=[-1,1]^2$ and the inside problem is $-\delta \Delta u + 2 \frac{\partial u}{\partial x} + \frac{\partial u}{\partial y} = f$. The source term $f$ and the Dirichelet boundary conditions are fixed to have an exact solution:
$u(x,y)=\left(1-e^{-\frac{1-x}{\delta}}\right)\left(1-e^{-\frac{1-y}{\delta}}\right)cos(\pi(x+y))$. We consider the two cases $\delta=0.1, 0.01$, and report in Figures \ref{Fig:boundarylayer}-\ref{Fig:boundarylayer_2} the results obtained.

    

\subsection{Anisotropic diffusion}
    \label{sec:Anisotropic_diffusion}
\begin{figure}[H]
        \centering
        \includegraphics[width=0.45\linewidth]{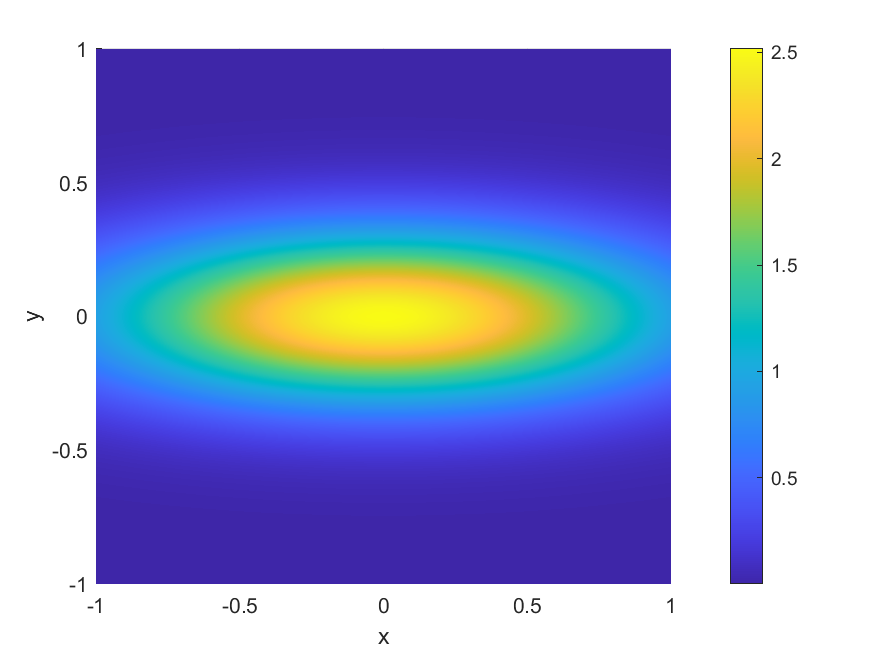}
        \includegraphics[width=0.45\linewidth]{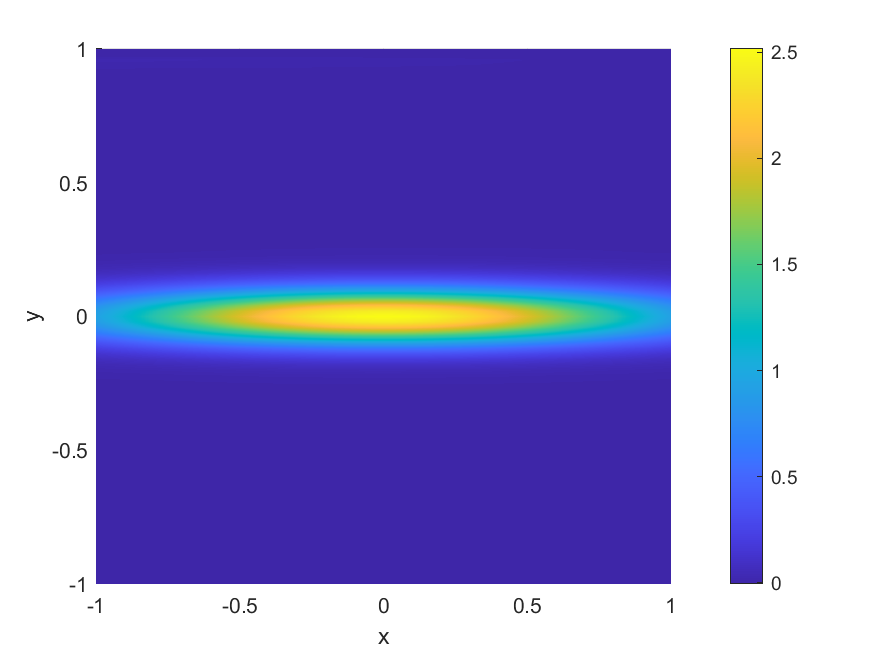}
        \includegraphics[width=0.45\linewidth]{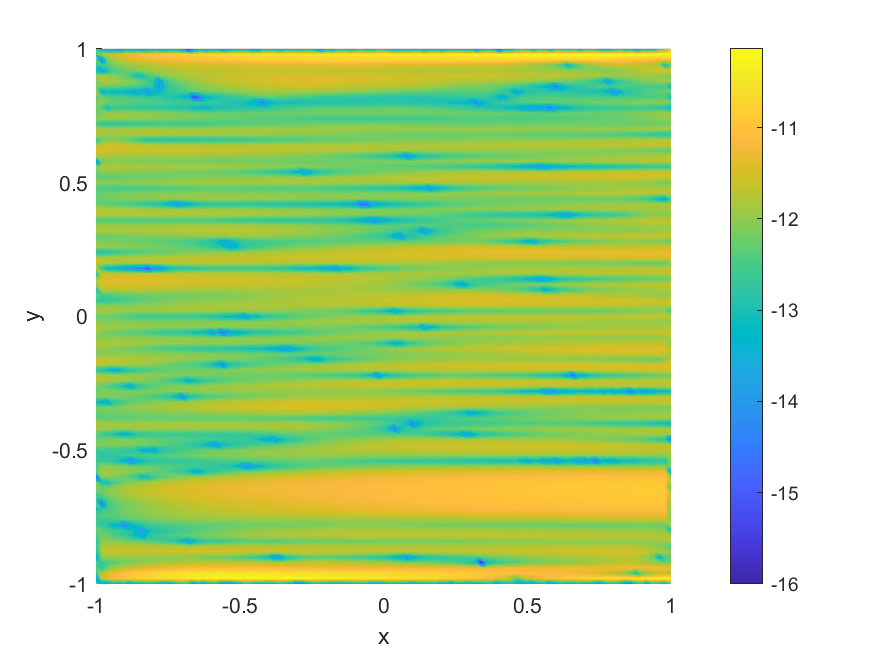}
        \includegraphics[width=0.45\linewidth]{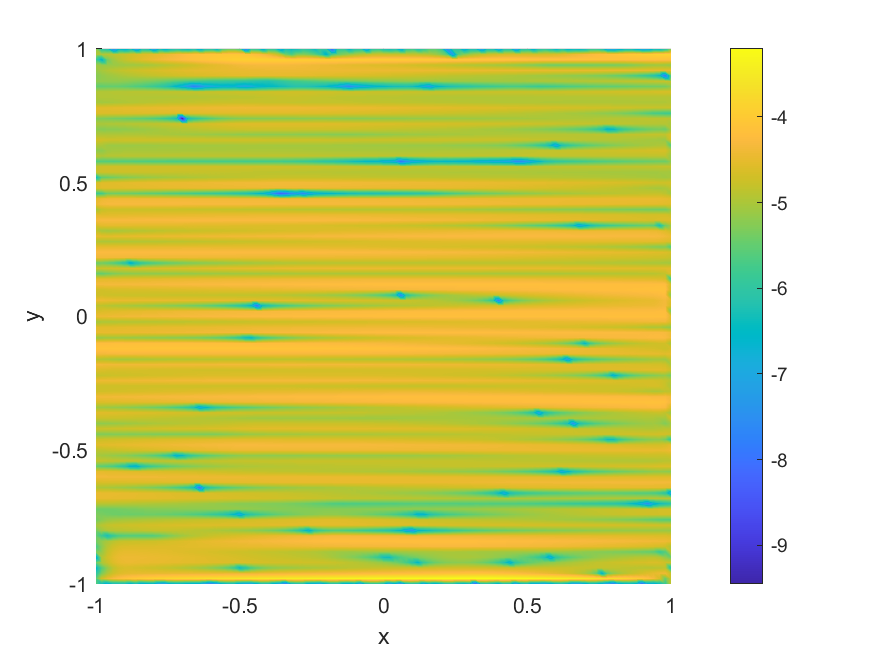}
        \caption{Anisotropic Diffusion problem. On the left, the case where $r = 10^4$, on the right $r=10^5$. The latter gives a stiff behavior. On the top, the exact solution, on the bottom computed errors in the case $N = 1000$.}
        \label{fig:anisotropic_diffusion}
    \end{figure}

\begin{figure}[H]
        \centering
        \includegraphics[width=0.48\textwidth]{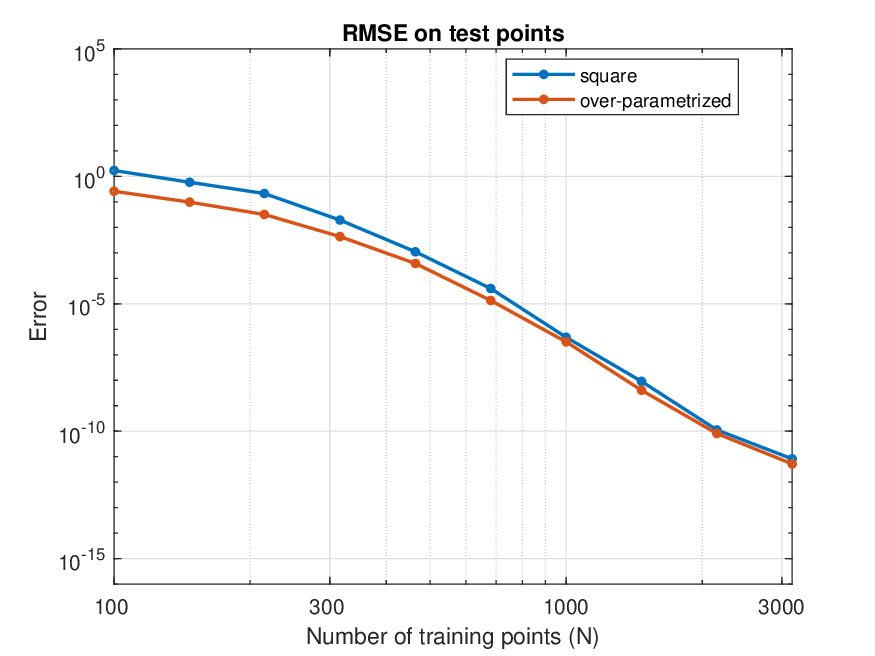}
        \includegraphics[width=0.48\textwidth]{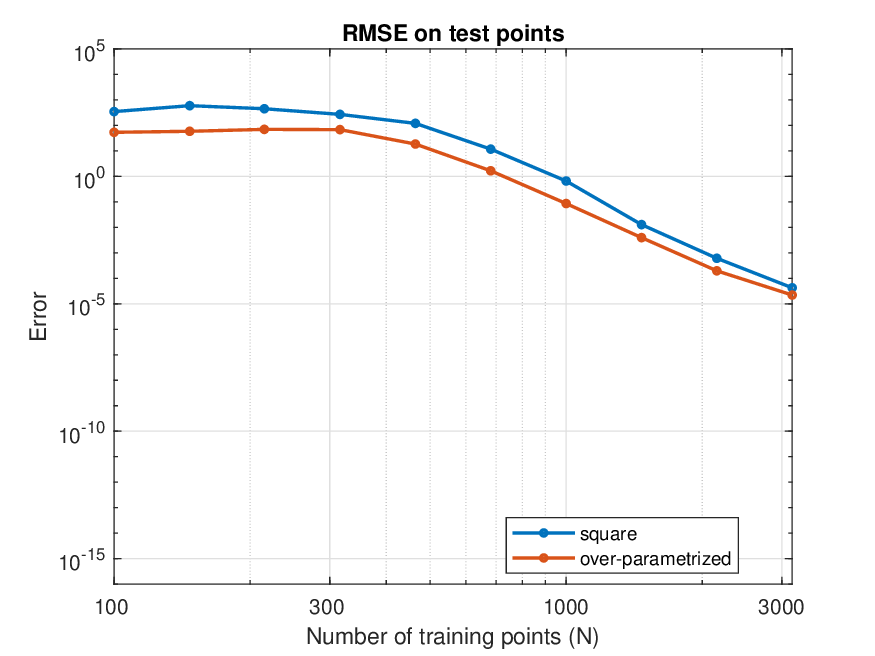}
        \caption{Convergence plots for the Anisotropic Diffusion problems. On the left, the case where $r=10^4$, on the right $r=10^5$. As it can be seen, the convergence strongly depends on the stiff behavior of the solution.}
        \label{fig:anisotropic_diffusion_2}

    \end{figure}

    In these examples, we consider diffusion problems with strong anisotropy. The domain is $\Omega=[-1,1]^2$ and the inside problem is $\frac{\partial^2 u}{\partial x^2} + \frac{1}{r}\frac{\partial^2 u}{\partial y^2} = f$. The source term $f$ and the Dirichlet boundary conditions are fixed to have an exact solution:
    $u(x,y)=\frac{e^{-(x^2 + \tau r y^2)}}{4\pi\sqrt{\tau}}$. We fix $\tau = 10^{-3}$ and consider the two cases $r=10^4, 10^5$.\\
    The parameter $r$ that appears both in the operator and in the exponent of the solution is used to specify the degree of anisotropy in the problem and the sharpness of the gradients in the solution.\\
    For $r=10^3$, the solution is a gaussian, and for increasing values of $r$ it gets tighter along the $y$ direction.
    We report in Figures \ref{fig:anisotropic_diffusion}-\ref{fig:anisotropic_diffusion_2} the results obtained.\\
    
\subsection{Reentrant corner}\label{sect:corner}
    \label{sec:Reentrant_corner}
\begin{figure}[H]
        \centering
        \includegraphics[width=0.45\linewidth]{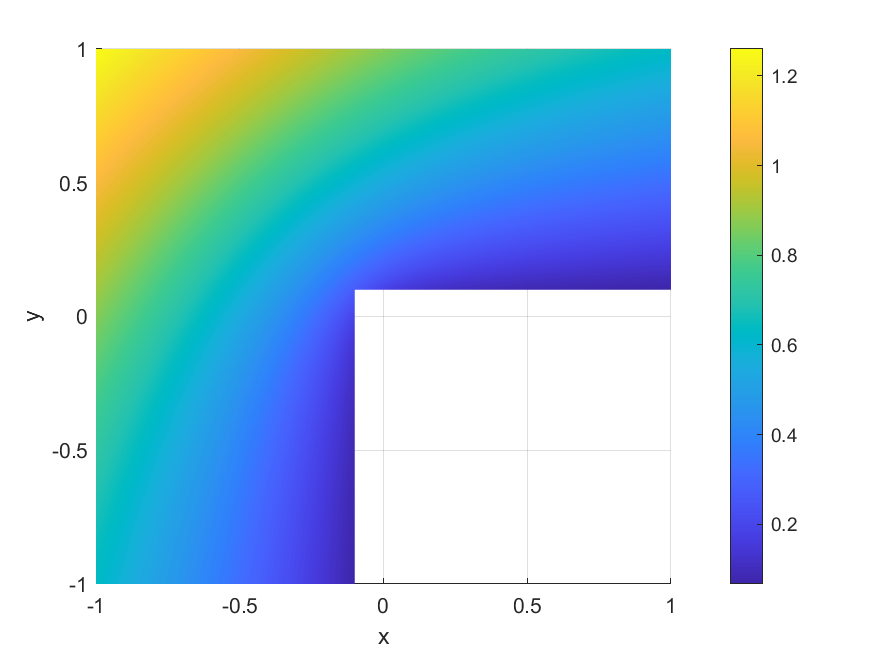}
        \includegraphics[width=0.45\linewidth]{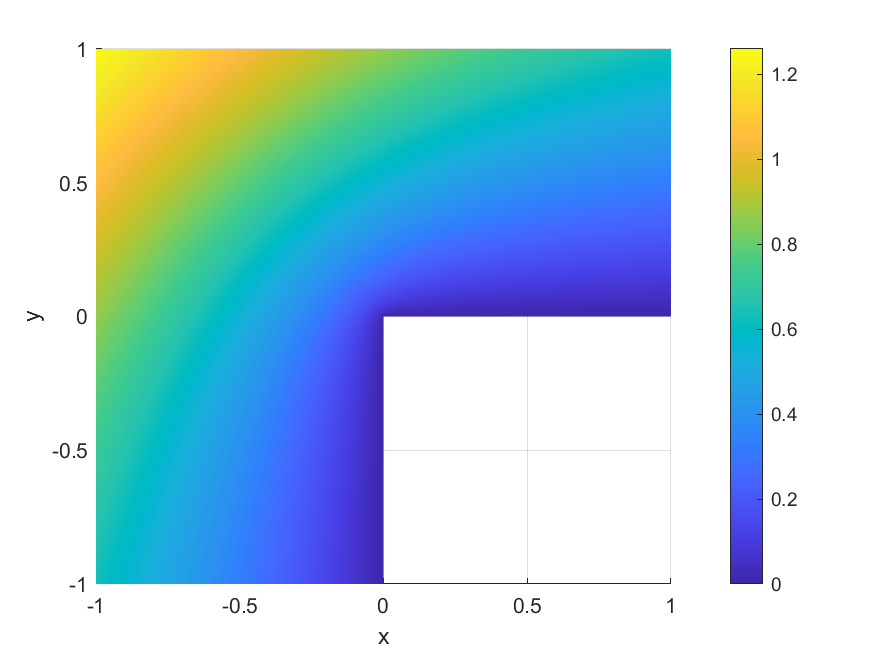}
        \includegraphics[width=0.45\linewidth]{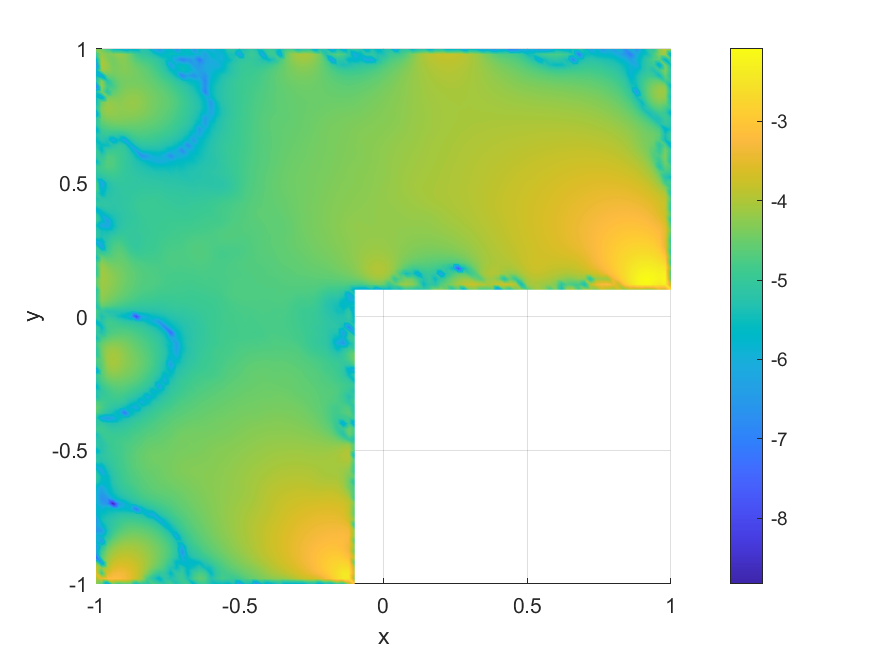}
        \includegraphics[width=0.45\linewidth]{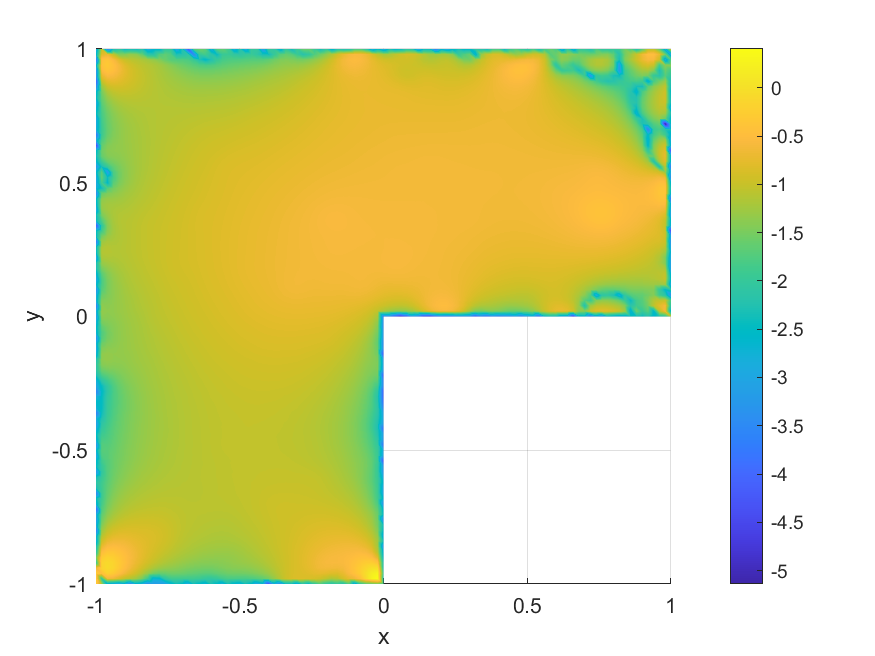}
        \caption{Reentrant Corner problem. On the left, the case where the singularity point is not in the domain, on the right, the case where the singularity is included. In the latter the gradient of the solution diverges in the origin. On the top, the exact solution, on the bottom computed errors in the case $N = 1000$.}
        \label{fig:reentrant_corner}
    \end{figure}

\begin{figure}[H]
        \centering
        \includegraphics[width=0.45\linewidth]{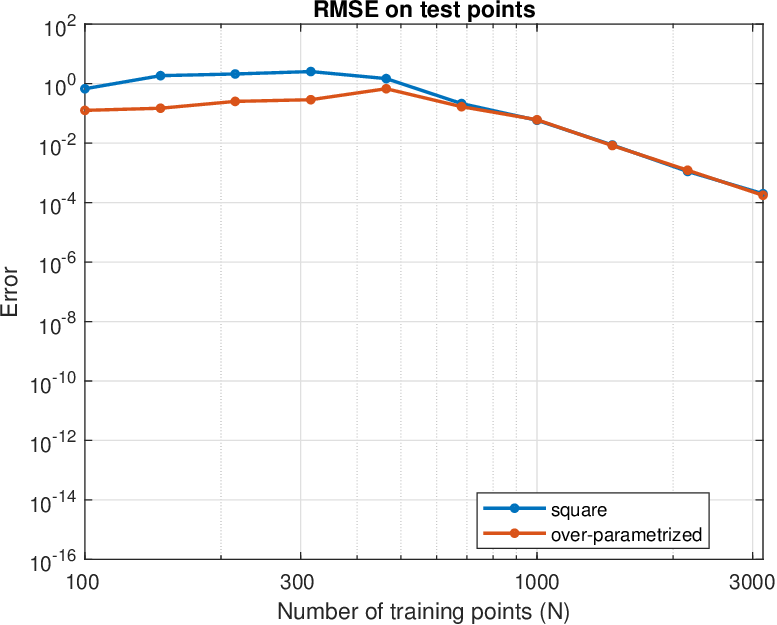}
        \includegraphics[width=0.45\linewidth]{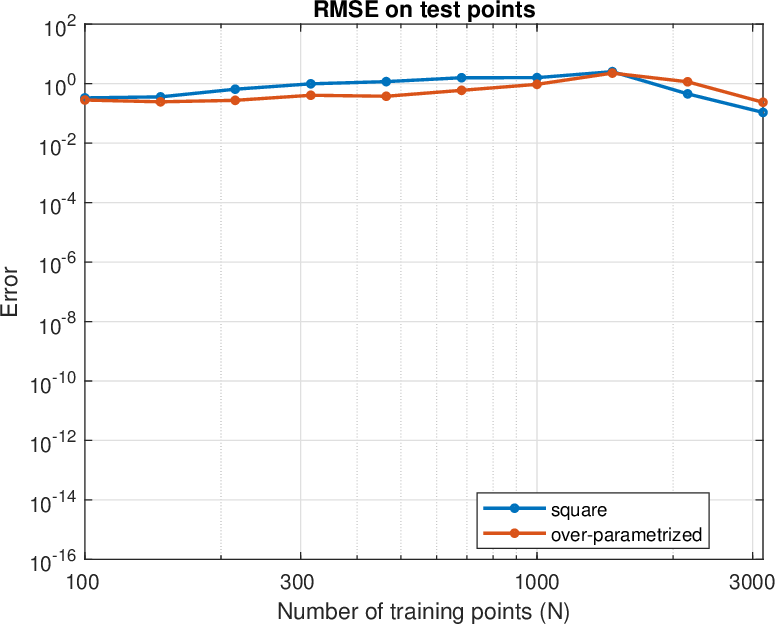}
        
        \caption{Convergence plots for the Reentrant Corner problems. On the left, the case without the singularity, on the right, the case with the singularity. As it can be seen, in the presence of the singularity the proposed method is not able to approximate the solution.}
        \label{fig:reentrant_corner_2}
    \end{figure}

    In these examples, we consider a Laplace problem whose solution has a singularity in the origin. The domain is $\Omega = [-1, 1]^2$ without a subregion on the bottom right and the inside problem is $\Delta u = 0$. The source term $f$ and the Dirichlet boundary conditions are fixed to have an exact solution (in polar coordinates): $r^\frac{2}{3} \, sin(\frac{2}{3} \theta)$. We consider these two cases: first we exclude the singularity point by removing a region that is slightly bigger than a quarter of the domain, and then we include it and this results in the well known L-shaped domain. In this case the solution has a singularity in the origin, therefore the magnitude of the Laplacian goes to infinity near the origin and the linear system that needs to be solved is extremely ill-conditioned.
    We report in Figures \ref{fig:reentrant_corner}-\ref{fig:reentrant_corner_2} the results obtained.

\subsection{1D Burgers-Fisher equation}
\begin{figure}[H]
    \centering
    \includegraphics[width=0.45\linewidth]{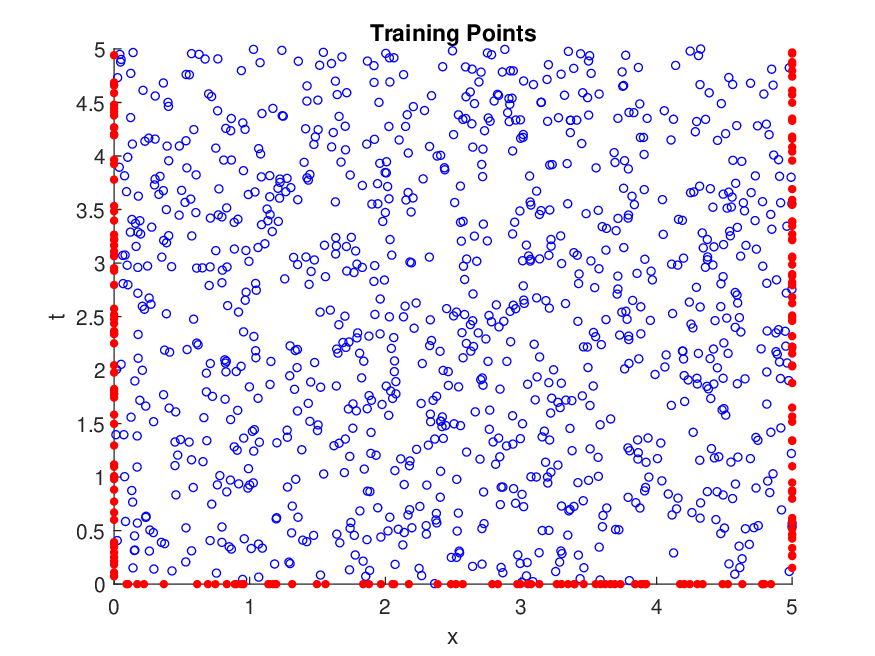}
    \includegraphics[width=0.45\linewidth]{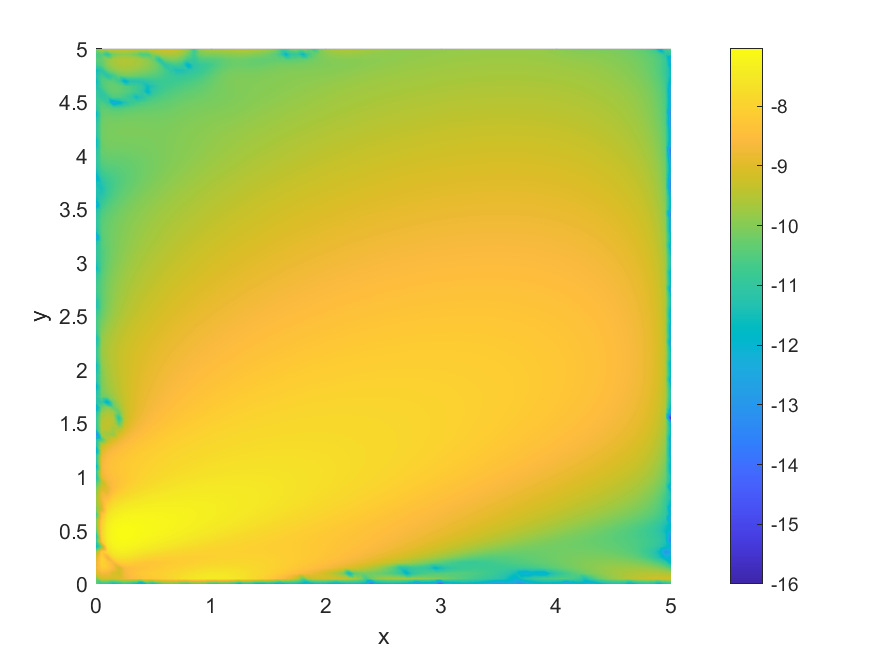}
    \caption{Example for the over-parametrized network with $N=1000$. On the left the distribution of training points, on the right the computed error on logarithmic scale}
    \label{fig:burgers_fisher_error}
\end{figure}

In this first time-dependent example, we focus on a problem with nonlinearities in the transport and reaction terms:
\begin{equation}
    \frac{\partial u}{\partial t} + u\frac{\partial u}{\partial x} = \frac{\partial^2 u}{\partial x^2} + u(1-u)\ .
\end{equation}
The spatial domain is $\Omega=[0,5]$ and the equation is solved for values of $t$ in the interval $[0,5]$, resulting in a square domain for a space-time collocation approach. The Dirichlet boundary conditions are fixed to have an exact solution: $u(t) = \frac{1}{2} + \frac{1}{2}tanh\left( \frac{5}{8}t - \frac{x}{4} \right)$ .\\
We report in Figure \ref{fig:burgers_fisher_error} the results obtained. It is important to note that the red dots on the "bottom" of the domain are points in which the initial condition is enforced, while on the two sides the boundary conditions are enforced as in the previous examples.\\
The nonlinearity has been solved by Gauss-Newton iterations, with a control on the norm of the update. This method converges in few iterations and kepps the accuracy that was obtained in the linear cases.

\begin{figure}[H]
    \centering
    \includegraphics[width=0.5\linewidth]{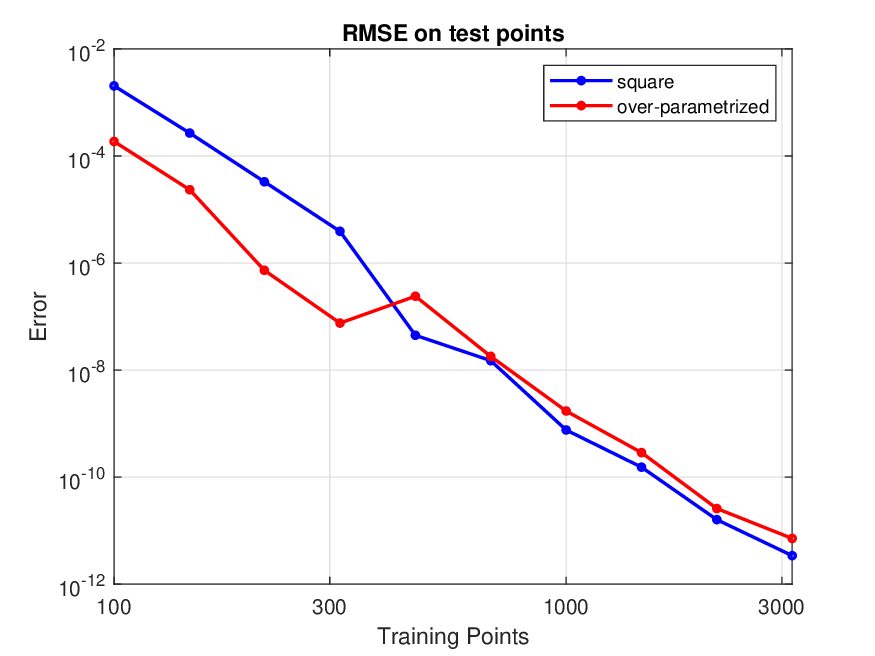}
    \caption{Convergence plot for the square (same number of neurons as training points) and over-parametrized cases. }
    \label{fig:burgers_fisher_convergence}
\end{figure}

Figure \ref{fig:burgers_fisher_convergence} shows the convergence of the error on test points as the number of training points and neurons in the hidden layer increase. Notably, the convergence rate for this nonlinear example is similar to other linear problems with similar regularity in the solution.

\subsection{2D Nonlinear Diffusion}
\begin{figure}[H]
    \centering
    \includegraphics[width=0.5\linewidth]{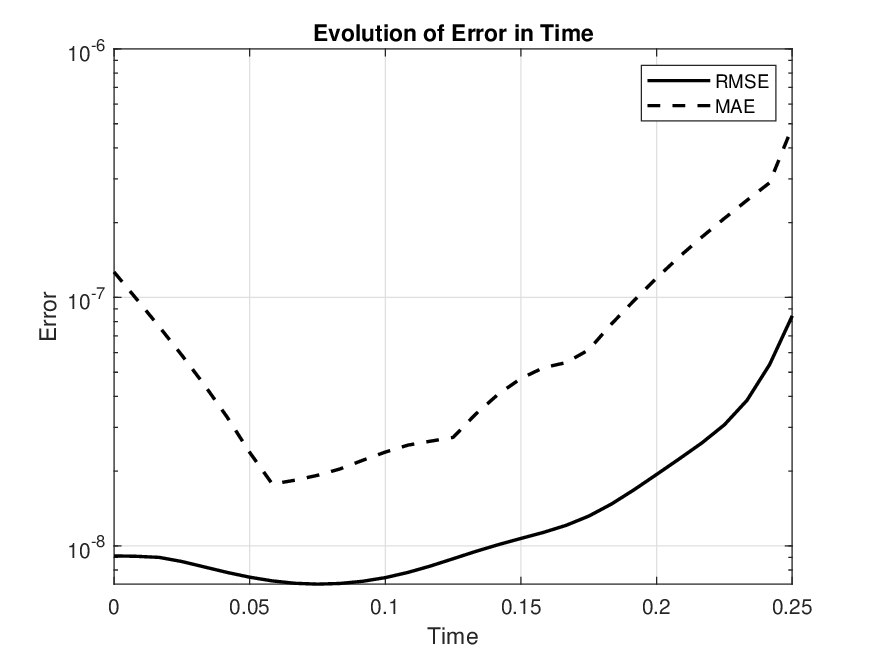}
    \caption{The error is computed on the spatial domain at each instant, then the values are plotted as a function of time}
    \label{fig:nonlinear_diffusion_error}
\end{figure}

\begin{figure}[H]
    \centering
    \includegraphics[width=0.5\linewidth]{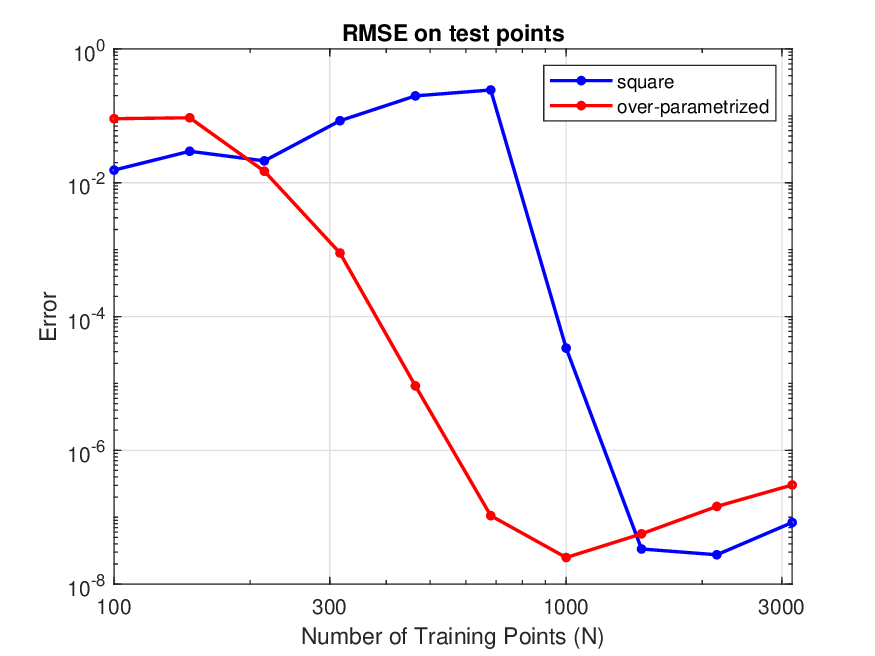}
    \caption{Convergence plot for the square and over-parametrized cases.}
    \label{fig:nonlinear_diffusion_convergence}
\end{figure}

In the last example, we consider a nonlinear diffusion:
$\frac{\partial u}{\partial t} = u\Delta u$
The spatial domain is $\Omega=[0,1]^2$ and the equation is solved for values of $t$ in the interval $[0,0.25]$. The Dirichlet boundary conditions are fixed to have an exact solution: $u = \frac{1}{1-t}\left( \frac{x^2+y^2+1}{4}+x\,y+x+y \right)$\\
We report in Figure \ref{fig:nonlinear_diffusion_error} the evolution in time of the spatial error obtained. A little increment as time goes on is expected since the gradient of the solution also increases\\
Figure \ref{fig:nonlinear_diffusion_convergence} shows the convergence of the error on test points as the number of training points and neurons in the hidden layer increase. The method is able to achieve good accuracy with a similar amount of training points, despite the fact that we are working on a 3D domain. From the figure, it can be noticed that the algorithm starts to converge when a sufficient number of points is taken and reaches a good accuracy rapidly, while it cannot gain higher accuracy when increasing the number of degrees of freedom.


\section{Conclusions}\label{sect:Conclusions}
In the present work, we have presented a new method for the numerical resolution of linear and nonlinear differential problems using single hidden-layer neural networks. This strategy is able to consider boundary conditions in an efficient way and mitigates some negative aspects of the previous strategies. The proposed LSE-ELM method is a global (space-time) collocation strategy, and this gives many pros and some cons. The pros include the possibility to consider points outside the domain in the resolution, so as to move the computation in the direction of physically based extrapolation. In fact, in the resolution of PDEs from the point of view of training, one can consider that the usual boundary conditions are good "supervised data", while the satisfaction of the differential problem can be a "physical constraint". Then, the resolution both inside the domain and outside can be an important feature in applications where extrapolation is required. The cons include the resolution of a full matrix and the loss of accuracy when the problem is less regular. The latter is an issue that can be fixed only with local strategies or the use of ad hoc activation functions. For the first contra, we have noticed with some computations that the spectra are better than the direct competitors, and we believe that some adaptation can further improve such a task. The constructions in the present paper are done in the most general case, and no adaptivity has been included, as an example with the use of Null Rules \cite{NL1,NL2} one could adapt the shape of the functions. However, the strategy has been used in such a setting in order to emphasize the ability of the method to deal with randomly chosen points and activation functions.



\bibliographystyle{plain}
\bibliography{references.bib}
\end{document}